\documentclass[a4paper]{amsart}
\usepackage[T1]{fontenc}

\usepackage{amssymb}

\usepackage{amsthm,mathtools}
\usepackage[active]{srcltx}
\usepackage{color}
\usepackage{enumerate}
\usepackage[all]{xy}
\usepackage{enumitem}

\usepackage[all]{xy}

\renewcommand{\k}{\Bbbk}

\newcommand{\nc}{\newcommand}

\newcommand{\fd}{finite-dimensional}

\newcommand{\ot}{\otimes}
\newcommand{\cf}{{\it cf.~}}
\newcommand{\eg}{{\it e.g.~}}
\newcommand{\ie}{{\it i.e., }}
\newcommand{\loc}{{\it loc.cit.~}}

\newcommand{\bs}{\boldsymbol}
\newcommand{\ydkg}{{}^{\k G}_{\k G}\mathcal{YD}}
\newcommand{\ydl}{{}^{\k \Lambda}_{\k \Lambda}\mathcal{YD}}
\newcommand{\ydkga}{{}^{\Gamma}_{\Gamma}\mathcal{YD}}
\newcommand{\ydga}{{}^{\k\Gamma}_{\k\Gamma}\mathcal{YD}}

\newcommand{\ydh}{{}^{H}_{H}\mathcal{YD}}
\newcommand{\ydL}{{}^{L}_{L}\mathcal{YD}}

\newcommand{\ydg}{{}^{G}_{G}\mathcal{YD}}
\newcommand{\ydgb}{{}^{\Gb}_{\Gb}\mathcal{YD}}

\newcommand{\ydk}{{}^{K}_{K}\mathcal{YD}}
\newcommand{\yda}{{}^{A}_{A}\mathcal{YD}}
\newcommand{\ydad}{{}^{A^*}_{A^*}\mathcal{YD}}
\newcommand{\ydhd}{{}^{H^*}_{H^*}\mathcal{YD}}
\newcommand{\ydkag}{{}^{\k^G}_{\k^G}\mathcal{YD}}
\newcommand{\fk}{\mathcal{FK}}

\newcommand{\tx}{\texttt{x}}

\newcommand{\rg}{\rangle}
\renewcommand{\lg}{\langle}
\renewcommand{\star}{*}

\newcommand{\Hom}{\mathrm{Hom}}
\newcommand{\one}{\mathbf{1}}
\newcommand{\id}{\operatorname{id}}
\newcommand{\FK}{\operatorname{FK}}
\newcommand{\ad}{\operatorname{ad}}
\newcommand{\sgn}{\operatorname{sign}}
\newcommand{\Ss}{\mathcal{S}}
\newcommand{\GL}{\operatorname{GL}}
\newcommand{\Alg}{\operatorname{Alg}}

\newcommand{\gr}{\operatorname{gr}}

\newcommand{\ev}{\operatorname{ev}}
\newcommand{\coev}{\operatorname{coev}}
\newcommand{\ord}{\operatorname{ord}}

\newcommand{\Cleft}{\operatorname{Cleft}}

\newcommand{\Aff}{\operatorname{Aff}}
\newcommand{\co}{\operatorname{co}}

\renewcommand{\mod}{\operatorname{-mod}}
\newcommand{\comod}{\operatorname{-comod}}
\newcommand{\svect}{\operatorname{SuperVect}}
\newcommand{\Sym}{\mathbb{S}}
\newcommand{\N}{\mathbb{N}}
\newcommand{\Z}{\mathbb{Z}}

\newcommand{\I}{\mathbb{I}}
\newcommand{\Q}{\mathcal{Q}}
\newcommand{\F}{\mathbb{F}}

\newcommand{\q}{\mathbf{q}}

\newcommand{\Gb}{\mathbf{G}}
\newcommand{\mI}{\mathcal{I}}

\newcommand{\C}{\mathcal{C}}
\renewcommand{\H}{\mathcal{H}}
\newcommand{\K}{\mathcal{K}}
\newcommand{\A}{\mathcal{A}}
\newcommand{\B}{\mathcal{B}}
\newcommand{\E}{\mathcal{E}}
\renewcommand{\O}{\mathcal{O}}
\newcommand{\D}{\mathcal{D}}
\newcommand{\R}{\mathcal{R}}
\newcommand{\J}{\mathcal{J}}
\newcommand{\Dd}{\mathbb{D}}
\renewcommand{\L}{\mathfrak{u}}
\newcommand{\Aa}{\mathbb{A}}
\newcommand{\mJ}{\mathrm{J}}

\renewcommand{\b}{\mathfrak{B}}

\newcommand{\eps}{\epsilon}

\def\pf{\begin{proof}}
\def\epf{\end{proof}}

\nc{\ben}{\begin{enumerate}[(i)]}
\nc{\een}{\end{enumerate}}

\def\trid{\vartriangleright}
\def\dirt{\vartriangleleft}
\numberwithin{equation}{section}
\theoremstyle{plain}
\newtheorem{theorem}{Theorem}[section]
\newtheorem{lemma}[theorem]{Lemma}

\newtheorem{proposition}[theorem]{Proposition}
\newtheorem{corollary}[theorem]{Corollary}

\theoremstyle{remark}
\newtheorem{remark}[theorem]{Remark}
\newtheorem*{acknowledgement*}{Acknowledgment}

\newtheorem{example}[theorem]{Example}
\newtheorem{definition}[theorem]{Definition}

\title{Twisting Hopf algebras from cocycle deformations}
\author[Andruskiewitsch; Garc\'ia Iglesias]{Nicol\'as Andruskiewitsch;
 Agust\'in Garc\'ia Iglesias}

\address{FaMAF-CIEM (CONICET), Universidad Nacional de C\'ordoba,
	Medina A\-llen\-de s/n, Ciudad Universitaria, 5000 C\' ordoba, Rep\'
	ublica Argentina.} \email{(andrus|aigarcia)@famaf.unc.edu.ar}

\thanks{The work was partially supported by CONICET,
	 Secyt (UNC), the MathAmSud project GR2HOPF}

\begin{document}

\begin{abstract}
Let $H$ be a Hopf algebra. Any finite-dimensional lifting of $V\in\ydh$ arising as a cocycle deformation of $A=\b(V)\#H$ defines a twist in the Hopf algebra $A^*$, via dualization. We follow this recipe to write down explicit examples and show that it extends known techniques for defining twists.
We also contribute with a detailed survey about twists in braided categories.
\end{abstract}

\maketitle

\section{Introduction}\label{sec:general}

Let $\k$ be an algebraically closed field of characteristic zero.
The classification of finite-dimensional Hopf algebras is an active area of research, see  \cite{A} and references therein. In particular, significant results on pointed Hopf algebras are available. The study of more general Hopf algebras is less developed, and new examples could shed more light on the subject. A specially fruitful approach is that of twisting: This procedure leaves the algebra structure unchanged and perturbs the comultiplication 
by conjugation by a so-called twist; originally due to Drinfeld
\cite{Dr}, it was adapted to Hopf algebras in \cite{R}. 

Given a Hopf algebra $H$, a twist is an invertible element $J\in H\ot H$ such that
\begin{align}\label{eqn:twist}
(1\ot J)(\id\ot\Delta)(J)=(J\ot 1) (\Delta\ot\id)(J), \quad (\id\ot\eps)(J)=(\eps\ot\id)(J)=1. 
\end{align}
Twists for braided Hopf algebras, namely {\it braided twists} \cite{Dr,R}, are defined analogously, see Definition \ref{def:braided-twists-cocycles} (i).
We present a detailed survey of some well-known facts abut twist in braided categories, which are hard to find in the literature. We obtain some new results that contribute to our goal of finding new examples of Hopf algebras via twisting. We show that 
if $R\in\ydh$ is a braided Hopf algebra, $\J\in R\ot R$ is a braided twist and  $J\in H\ot H$ is a twist, then 
\begin{align}\label{eqn:def:twist-bosonization-intro}
\J\# J:=J^1_{(1)}\cdot \J^1\# J^1_{(2)} \J^2_{(-1)}  \ot J^2_{(1)}\cdot \J^2_{(0)} \# J^2_{(2)}
\end{align}
is a twist for $R\#H$, see Proposition \ref{pro:JF}. We show that many existing examples can be described in this way.

\smallbreak
Twists for group algebras $\k G$, $G$ a finite group, are classified in \cite{Mv}, \cf also \cite{EG}. 
See \cite{AEGN} for the analogous result in positive characteristic.

\smallbreak
Twists for supergroup algebras, \ie Hopf algebras $R=\k G\ltimes \Lambda(V)\in\svect$, were classified in \cite{AEG}: they arise as $\J=\exp (r/2)$, $r\in S^2(V)$. This  is equivalent to the classification of the twists for triangular Hopf algebras of {\it rank} $\leq 2$, which arise as $\J\#1$ \cf \cite[Corollary 3.3.3 \& Theorem 5.1.1]{AEG}.

\smallbreak
Examples of twists for Hopf algebras of the form $\B(V)\#\k G$, where $V$ is a  quantum linear space over the finite group $G$ with $\dim \B(V) < \infty$, were presented in \cite{Mo}, extending \cite{AEG}. They arise as $\J_{\D}\# J$, $J$ a twist for $H$ and a braided twist $\J_D$ depending on a certain family of scalars $\D$. We review these examples in \S \ref{sec:qls}.

\smallbreak

In this article we consider any $V\in\ydg$ with $\dim \B(V) < \infty$ and construct twists for Hopf algebras of the form  $\B(V)\# \k G$ or $\B(V)\# \k^G$, extending \cite{Mo}. 

\smallbreak
Our approach is as follows:

\begin{enumerate} [leftmargin=*]\renewcommand{\theenumi}{\alph{enumi}}\renewcommand{\labelenumi}{(\theenumi)}
	\item If $\dim H <\infty$, then an invertible $J\in H\ot H$ is a twist if and only if the map $\sigma: H^*\ot H^*\to \k$,  
	$f\ot g\mapsto (f\ot g)(J)$ is a Hopf 2-cocycle, \cf \S \ref{sec:cocycle}.

\smallbreak
\item Let $K$ be a \fd{} cosemisimple Hopf algebra, $V\in\ydk$ such that $\dim \B(V) < \infty$ and set  $\K=\b(V)\# K$.
Let $\L$ be a {\it lifting} of $V$, \ie $\gr \L \simeq \K$,
where $\gr \L$ is the graded Hopf algebra corresponding to the coradical filtration. 
In all known cases, $\L$ is a cocycle deformation of $\K$, that is 
\begin{align}\label{eq:cocycle-def}
\L & = \K_\sigma, \qquad  \sigma \text{ a Hopf 2-cocycle.}
\end{align} 
Then the transpose ${}^t\sigma:\k\to \K^*\ot \K^*$ of $\sigma$ defines a twist $J(\sigma)={}^t\sigma(1)$ for $\A = \K^* \simeq \b(V^*)\# K^*$.

\smallbreak
\item\label{c} Let $L$ be a \fd{} cosemisimple Hopf algebra and $U\in\ydL$ such that $\dim \B(U) < \infty$. We denote by $\cdot :L\ot U\to U$, resp. $\delta:U\to L\ot U$ the action, resp. the coaction. Assume that
\begin{align}\label{eq:braided-subspace}
\begin{split}
A\subseteq L & \text{ is a Hopf subalgebra}, \\ 
W\subseteq U & \text{ is an $A$-stable, $A$-co-stable braided subspace},
\end{split}
\end{align}
that is, $A\cdot W\subseteq W$ and $\delta(W)\subseteq A\ot W$, so that  $W \in \yda$.

\smallbreak
Let $K = A^*$, $V = W^*$. If \eqref{eq:cocycle-def} holds, then $J(\sigma)$ is a twist for $\H= \b(U)\# L$.
\end{enumerate}

In other words, we construct twists by dualizing Hopf 2-cocycles, obtained indirectly from theoretical arguments 
in lifting theory. This is indeed a fertile source of Hopf 2-cocycles:

\begin{enumerate} [leftmargin=*]\renewcommand{\theenumi}{\roman{enumi}}\renewcommand{\labelenumi}{(\theenumi)}
	
	\smallbreak
	\item Let $\Gamma$ be a finite abelian group whose order is not divisible by 2,3,5,7. If $\L$ is a \fd{} pointed 
	Hopf algebra with $G(\L) \simeq \Gamma$, then there exists $V \in \ydkga$ such that $\L$ is a lifting of $V$; all such liftings are classified	\cite{AS} and all of them are cocycle deformations \cite{M2}.
	
	\smallbreak
	\item  Let $K$ be a cosemisimple Hopf algebra and $V$ a braided vector space of Cartan type $A_n$, with a principal realization in $\ydk$. Then all liftings of $V$ are classified	and all of them are cocycle deformations \cite{AAG}. This last statement holds true for all $V\in \ydk$ of diagonal type \cite{AGI}.
	
	\smallbreak
	\item Let $G$ be either $\Sym_3$, $\Sym_4$ or $\mathbb{D}_{4t}$. If $\L$ is a \fd{} pointed 
	Hopf algebra with $G(\L) \simeq G$, then there exists $V \in \ydg$ such that $\L$ is a lifting of $V$. All such liftings are classified	\cite{AG,GG,FG} and all of them are cocycle deformations \cite{GMo,FG}. This also holds for certain semidirect products $\Z_p\rtimes \Z_q$ \cite{GV}. 
\end{enumerate}

\smallbreak \label{page:intro}A setting where \eqref{c} above holds is, for instance, when $U$ is a braided vector space of diagonal type, with a principal realization $U\in\ydL$, so there is $\Gamma\leq Z(G(H))$ such that $W = U\in\yda$, $A = \k\Gamma$. If $U$ is of diagonal type, then \eqref{eq:cocycle-def} holds.

\smallbreak 
The paper is organized as follows: In \S \ref{sec:prelim} we explain the notation we shall use throughout the text. In the survey Section \ref{sec:twist-and-cocycles} we 
recall the notions of a twist and Hopf 2-cocycle on a Hopf algebra and their braided analogues. In the braided setting, in 
particular, we compile in several lemmata some well-known facts about braided twists and Hopf 2-cocycles which, however, are hard to 
find in the literature. We relate these facts to the Hopf 2-cocycles arising from liftings Hopf algebras in \S \ref{sec:twist-and-lift} and explain how to produce twists out of a given deformation. In  \S \ref{sec:main} we specialize this procedure on examples from the theory of pointed and copointed Hopf algebras. 

\section*{Acknowledgments}
We thank the referee for a careful reading of our article. His/Her suggestions helped us to improve the exposition.

\section{Preliminaries}\label{sec:prelim}

We set $\I_n:=\{1,\dots, n\}$, $n\in\N$, and omit the subscript when it is clear from the context. We denote by $\Sym_n$, resp. $\Aa_n$, the symmetric, resp. alternating, group on $n$ letters. Also, $\Dd_n$ denotes the dihedral group of order $2n$. If $V$ is a $\k$-vector space, then $V^{\ot n}$ denotes the iterated tensor product of $V$ with itself, $n$ times. We denote the evaluation by $\lg \,,\,\rg_V:V^*\ot V\to\k$, dropping the subscript when it is clear from the context. If $V,W$ are finite-dimensional $\k$-vector spaces, then we shall use the (inverse of the) identification $W^*\ot V^*\stackrel{\sim}{\longrightarrow} (V\ot W)^*$ given by 
\begin{align}\label{eqn:dual-identification}
g\ot f\longmapsto \left(v\ot w\mapsto \lg f,v\rg_V\lg g,w\rg_W, \quad v\in V,w\in W\right), \quad f\in V^*, g\in W^*.
\end{align}
If $G$ is a group, then $Z(G)$ denotes the center of $G$ and $\widehat{G}$ the group of characters (1-dimensional representations) of $G$. For a $\k$-algebra $A$, we write $A^{\times}$ for the group of its invertible elements. If $C$ is a
coalgebra, then we denote by $f\star g \in \Hom(C,A)$ the {\it convolution product} of $f,g\in\Hom(C,A)$.

Let $H$ be a Hopf algebra, with multiplication  $m$, comultiplication $\Delta$ and antipode $\Ss$ (always assumed bijective).
We use Sweedler's notation for comultiplication and coaction. 
We denote by $\Alg(H,\k)$ the group of 1-dimensional representations of $H$. We shall denote by $H\mod$, resp. $H\comod$ the categories of left $H$-modules, resp. left $H$-comodules. Recall that $H$ is said to have the {\it Chevalley property} \cite{AEG} when the tensor product of any two simple $H$-modules is semisimple. If $M\in H\mod$ and $\chi\in\Alg(H,\k)$, then $M^{\chi}$ denotes the isotypic component of 
type $\chi$. We denote by $G(H)$ the group of group-like elements of $H$. If $V\in H\comod$ and $g\in G(H)$, then we set $V_g=\{x\in V:x_{(-1)}\ot x_{(0)}=g\ot x\}$.

We write $H_0$ for the coradical of $H$; $(H_n)_{n\geq 0}$ denotes the coradical filtration and $\gr H=\oplus_{i\geq 0} H_i/H_{i-1}$ is the associated graded coalgebra; $H_{-1}:=0$. If $H_0$ is a Hopf subalgebra, then $\gr H$ is a graded Hopf algebra. 
Recall that $H$ is {\it pointed} when  $H_0\simeq \k G(H)$ and {\it copointed} when $H_0$ is isomorphic to the algebra of functions 
$\k^G$ of a non-abelian group $G$ \cite{AV}. 

We  denote by $\ydh$ the category of Yetter-Drinfeld modules over $H$. If $V\in\ydh$, then the {\it Nichols 
	algebra} $\b(V)$ is a graded Hopf algebra in $\ydh$ generated in degree one by $V$, 
	which is the subspace of its primitive elements. A lifting of $V$ is a Hopf algebra $L$ such that $\gr L\simeq \b(V)\# H$.
	See \cite{AS} for unexplained notation.

If $\dim H<\infty$, then $\ydh$ is rigid: if 
$V\in\ydh$, then $V^*\in\ydh$ with 
\begin{align}\label{eqn:dual-1}
 f_{(-1)} \lg f_{(0)},v \rg&:=\Ss^{-1}(v_{(-1)})\lg f,v_{(0)}\rg,  &
\lg h\cdot f,v\rg &:= \lg\ f,\Ss(h)\cdot v\rg.
\end{align}
for $h\in H$, $f\in V^*$, $v\in V$. If $V$ is \fd,
then also $V^*\in\ydhd$ via
\begin{align}\label{eqn:dual}
 \lg f_{[-1]},h\rg \lg f_{[0]},v \rg&:=\lg f,h\cdot v\rg,  &
\lg\alpha\rightharpoonup f,v\rg &:= \lg\alpha,v_{(-1)}\rg\lg\ f,v_{(0)}\rg.
\end{align}
for $\alpha\in H^*$, $h\in H$, $f\in V^*$, $v\in V$. 
As a result, if $V\in\ydh$, then $V\in\ydhd$ via
\begin{align}\label{eqn:H-to-H*}
\alpha\rightharpoonup v&=\lg\alpha,\Ss(v_{(-1)})\rg v_{(0)}, &\lg v_{[-1])},h\rg v_{[0]}&=\Ss^{-1}(h)\cdot v.
\end{align}
%

\subsection{Braided vector spaces and realizations} \label{sec:realizations}

Recall that a braided vector space  is a pair $(V,c)$, where $V$ is a vector space and $c\in\GL(V^{\ot 2})$ is a braiding, \ie
satisfies the {\it braid equation}
$(c\ot \id)(\id\ot c)(c\ot \id)=(\id\ot c)(c\ot \id)(\id\ot c)$.
The category $\ydh$ is braided; hence any $V\in\ydh$ becomes a  braided vector space. 
A \emph{realization} of $(V,c)$  over $H$ is a structure of 
Yetter-Drinfeld module on $V$ such that the braiding $c\in\GL(V^{\ot2})$ coincides with the braiding from $\ydh$ \cite{AG}.

\subsubsection{YD-pairs}\label{sec:ydpair}
A \emph{Yetter-Drinfeld pair} is  $(g, \chi)\in G(H)\times\Alg(H,\k)$ such that
$\chi(h)\,g = \chi(h_{(2)}) h_{(1)}\, g\, \Ss(h_{(3)})$, $h\in H$.
Hence $g\in Z(G(H))$ and $\k_g^{\chi}:=\k x\in\ydh$ with coaction
$x\mapsto  g\ot x$ and action $h\cdot x=\chi(h)x$, $h\in H$. 
Indeed $(g,\chi)$ gives a realization in $\ydh$ of the braided vector space $\k x$ with braiding $c=\chi(g)$.
Conversely, any one-dimensional Yetter-Drinfeld module over $H$ arises in this way.

\subsubsection{Diagonal type}\label{sec:diagonal-0}
A braided vector space  $(V,c)$ is of diagonal type if there is a basis 
$(x_i)_{i \in \I_\theta}$ and a {\it braiding matrix} $\q=(q_{ij})\in (\k^{\times})^{\theta\times\theta}$ such 
that  $c(x_i\ot x_j)=q_{ij}x_j\ot x_i$, $i,j\in \I_\theta$.
A \emph{principal realization} of $(V, c)$ is a collection $(\chi_j,g_j)_{j\in \I_\theta}$ of YD-pairs with 
$\chi_i(g_j)=q_{ji}$, $i,j\in \I_\theta$. Hence $V$ is realized in $\ydh$ and $\ydkga$, where
\begin{align}\label{eqn:support}
\Gamma=\lg g_1,\dots,g_\theta \rg\leq Z(G(H)).
\end{align}

\begin{remark}\label{rem:dual} 
	Set $\Lambda=\widehat\Gamma$ and let $W:=V^*\in\ydl$. Then $W$ is the braided vector 
	space associated to the transpose matrix $\,^t\q$. 
\end{remark}

\subsubsection{Rack type}\label{sec:rack}
Let $(X,\trid)$ be a finite rack and $q:X\times X\to \k$ be a rack 2-cocycle, that is a map $(i,j)\mapsto q_{ij}$ satisfying
\[
q_{i,j\rhd k}q_{j,k}=q_{i\rhd j,i\rhd k}q_{i,k}, \qquad i,j,k\in X.
\]
Then $(V, c^q)$ is a braided vector space, where $V$ has a  basis $(x_i)_{i\in X}$ and
$$
c^q(x_i\ot x_j)=q_{ij}\,x_{i\trid j}\ot x_i, \qquad i,j\in X.
$$
We refer to \cite{AG} for details. 
A {\it principal Yetter-Drinfeld realization} of $(X,q)$ over a finite group $G$, \cf \cite[Definition 3.2]{AG} is a collection
$(\cdot, g, (\chi_i)_{i\in X})$ where 
\begin{itemize}[leftmargin=*]\renewcommand{\labelitemi}{$\circ$}
 \item $\cdot$ is an action of $G$ on $X$;
\item $g:X\to G$ is a function such that $g_{h\cdot i} = hg_ {i}h^{-1}$
and $g_{i}\cdot j=i\rhd j$;
\item the family $(\chi_i)_{i\in X}$, with
$\chi_i:G\to\k^*$ satisfies
$\chi_i(ht)=\chi_i(t)\chi_{t\cdot i}(h)$, for all $i\in X$, $h,t\in
G$, together with $\chi_i(g_{j})=q_{ji}$, $i,j\in X$.
\end{itemize}
These data determine a  realization of $V$ as an object in $\ydkg$. 

\begin{example}
The affine rack  $\Q_{q, \omega}$, where $q$ is a prime power and  $\omega \in \F_q^{\times}$, is the set $\F_q$ 
with the rack multiplication $i\trid j:=(1-\omega)i+\omega j$, $\ i,j\in\F_q$. 
We denote by $\Aff(\F_b,\omega)$ the braided vector space associated to  $(\Q_{q, \omega}, -1)$.
\end{example}

\section{Twists and Hopf 2-cocycles}\label{sec:twist-and-cocycles}

We recall some well-known facts about twists and Hopf 2-cocycles. We include poofs of facts largely known in the folklore of the subject for which we could not find a proof. We fix a Hopf algebra $H$.

\subsection{Twists}

An  element $J\in (H\ot H)^\times$ is a {\it twist} if and only if it satisfies \eqref{eqn:twist} \cite{Dr,R}. If $J\in H\ot H$ is a twist, then $H^J$, the algebra $H$ with the comultiplication 
\begin{align}\label{eqn:twist-comult}
\Delta^J:=J\Delta J^{-1}, 
\end{align}
becomes a Hopf algebra.
Observe that if $J$ is a twist for $H$ and $H$ is a Hopf subalgebra of a Hopf algebra $K$, then $J$ is a twist for $K$. Also, if $J$ is a twist for $H$, then $J^{-1}$ is a twist for $H^J$ and $(H^J)^{J^{-1}}\simeq H$.
We shall use a Sweedler-type notation $J=J^1\ot J^2$. Also,  $J^{-1}:=J^{-1}\ot J^{-2}$.

When $H=\k G$, $G$ a finite group, twists in are classified up to gauge equivalence, \cf \cite{Mv, EG}, 
by conjugacy classes of pairs $(S,\alpha)$, where $S < G$ and $\alpha \in H^2(S,\k^\times)$ is non-degenerate. 
When $\Gamma$ is an abelian group, the twist $J$ corresponding to the (class)
$\alpha\in H^2(\Gamma,\k^\times)$ is given by
\begin{align}\label{eqn:twist-abelian}
J_\alpha =\sum_{\chi,\tau\in\widehat{\Gamma}}\alpha(\chi,\tau) e_\chi \ot e_\tau.
\end{align}
Here, $(e_\chi)_{\chi\in \widehat{\Gamma}}$ is a basis of 
orthogonal central idempotents of $\k\Gamma$.

\begin{remark}
We follow \cite[10.2.3, 10.2.4]{KS} to give the definition of twist and the twisted comultiplication $\Delta^J$. This definition differs (though it is equivalent to) the one considered in \cite{AEG,Mo}. We prefer ours as it is the one dualizing Hopf 2-cocycles, see \S \ref {sec:cocycle} next.
In particular, see Remark \ref{rem:duality}.
\end{remark}

The procedure of twisting has the following categorical interpretation given in \cite{S};
it holds
more generally for quasi-bialgebras \cite{EG1}. 
\begin{theorem}\label{thm:twist-versus-tensor} \cite[Theorem 6.1]{EG1}
Two finite-dimensional Hopf algebras $H$ and $H'$ are twist equivalent
if and only if  $H\mod$ and $H'\mod$ are tensor equivalent. \qed
\end{theorem}

\subsection{Hopf 2-cocycles}\label{sec:cocycle}
A Hopf 2-cocycle on $H$ (with values in $\k$) \cite{D2} is a 
convolution invertible linear map $\sigma:H\ot H\to \k$ satisfying
\begin{align}
\label{eqn:cocycle}
\sigma(x_{(1)},y_{(1)}) \sigma(x_{(2)}y_{(2)},z) & =\sigma(y_{(1)},z_{(1)}) \sigma(x,y_{(2)}z_{(2)}), \qquad x,y,z\in H.
\end{align}
and normalized by setting $\sigma(x,1)=\sigma(1,x)=\eps(x)$, $x\in H$. We denote by $Z^2(H)$ the set of normalized Hopf 2-cocycles on $H$. If $\sigma \in Z^2(H)$, then $H_\sigma$, the coalgebra $H$ with multiplication $m_\sigma =\sigma \star m\star\sigma^{-1}$, 
is again a Hopf algebra. We refer to $H_\sigma$ as a {\it cocycle deformation} of $H$. 

\begin{remark}
The Hopf algebra $H_\sigma$ was originally denoted $H^\sigma$ in \cite{D2} and the map $\sigma$ is called a cocycle, also a bialgebra cocycle. This is justified by the fact that when $H$ is cocommutative, then \eqref{eqn:cocycle} states that $\sigma$ is a 2-cocycle relative to $\k$ in the sense of \cite{Sw}. We point out that in \loc it is shown that the normalizing condition is equivalent to setting $\sigma(1,1)=1$.

When the context is clear, we shall refer to Hopf 2-cocycles as 2-cocycles or simply cocycles, indistinctly. These objects have been also called multiplicative cocycles in the literature. The corresponding Hopf algebra $H_\sigma$ has also been called a deformation, or a cocycle twist.
\end{remark}

\begin{remark}\label{rem:duality}
Let $H$ be a finite-dimensional Hopf algebra and consider the identification $H\ot H\simeq (H^*\ot H^*)^*$, via \eqref{eqn:dual-identification}.
Then $J\in H\ot H$ is a twist if and only if $f\ot g\mapsto \lg f\ot g,J\rg$, $f,g\in H^*$, defines a Hopf 2-cocycle on $H^*$.
\end{remark}

If $H=R\# K$, $R\in\ydk$, then $Z^2_K(H)$ will denote the set of $K$-balanced Hopf 2-cocycles, that is the set of those $\sigma\in 
Z^2(H)$ that factor through $H\ot_{K}H\to \k$ as $K$-bimodule homomorphisms. If $G\leq G(H)$ is a subgroup, then $Z^2_{G}(H):=Z^2_{\k G}(H)$.

\subsection{Cleft objects}\label{sec:cleft}
Let $A$ be a (right) $H$-comodule algebra, $B=A^{\co H}$ be the subalgebra of coinvariants. 
Recall that $A$ is a cleft extension of $B$ if there is an $H$-colinear convolution-invertible map $\gamma:H\to A$ \cf \cite{DT}. Such $\gamma$ is called a {\it cleaving map} and can be chosen so that $\gamma(1)=1$, in which case it is said to be a {\it section} \cite{Mont}. A cleft extension with trivial coinvariants $B=\k$ is called a {\it cleft object}. We denote by $\Cleft(H)$ the set of (isomorphism classes) of cleft objects for $H$.

If $\sigma:H\ot H\to\k$ is a Hopf 2-cocycle, then $m_{(\sigma)}:=\sigma\star m$ is as associative multiplication on the vector space $H$, which thus becomes a cleft object for $H$, denoted $H_{(\sigma)}$, with coaction induced by the comultiplication. Conversely, if $A$ is a cleft object with section $\gamma:H\to A$, then it determines a Hopf 2-cocycle $\sigma$ for $H$ by setting, \cf \cite[Theorem 11]{DT},
\begin{align}\label{eqn:sigma-gamma}
\sigma(x,y)=\gamma(x_{(1)})\gamma(y_{(1)})\gamma^{-1}(x_{(2)}y_{(2)}), \quad x,y\in H.
\end{align}

\subsection{Braided twists and braided Hopf 2-cocycles}\label{sec:braided-tc}
Let $(\C,\ot,\one)$ be a braided tensor category, see \cite{EGNO}, with braiding $(c_{V,W})_{V,W\in\C}$. 
Let $(A,m_A,u_A)$  be an algebra and let  $(C,\Delta_C,\eps_C)$ be a
coalgebra in $\C$. If $f,g\in\Hom(C,A)$, then we denote by
$f\star g \in \Hom(C,A)$ the {\it convolution product}, \ie the map given by the
commutativity of the diagram:
\begin{equation*}
\xymatrix{C \ar[rr]^{f\star g} \ar[d]_{\Delta_C} && A  \\
	C \ot C  \ar[rr]^{f\ot g}  && A \ot A. \ar[u]_{m_{A}}}
\end{equation*}
A map $f\in \Hom(C, A)$ is \emph{convolution invertible} iff there exists
$g\in \Hom(C,A)$ such that  $f \star g = g \star f=u_A\eps_C$.
Let also $X \in \C$, $f\in \Hom(\one,A)$, $h\in \Hom(X,A)$. We define $f\rightharpoonup h$, $h\leftharpoonup f \in \Hom(X,A)$  by the commutativity of the following diagrams:
\begin{align*}
&\xymatrix{X \ar[rr]^{f\rightharpoonup h} \ar[d]_{\simeq} && A  \\
	\one \ot X  \ar[rr]^{f\ot h}  && A \ot A, \ar[u]_{m_{A}}} &
&\xymatrix{X \ar[rr]^{h\leftharpoonup f} \ar[d]_{\simeq} && A  \\
	X\ot\one   \ar[rr]^{h\ot f}  && A \ot A. \ar[u]_{m_{A}}}
\end{align*}
Observe that if  $ X = \one$, then  $\leftharpoonup =  \rightharpoonup = \star$. 
If $(A,m_A)$, $(B,m_B)$ are algebras in $\C$, then we denote by $A\underline{\ot} B:=(A\ot B,
m_{A\ot B})$
the algebra in $\C$ with multiplication
\begin{equation}\label{eq:braided-mult}
m_{A\ot B}=(m_A\ot m_B)\circ(\id \ot \, c_{A,B} \ot \id).
\end{equation}
Given two coalgebras $(C,\Delta_C)$, $(D,\Delta_D)$ in $\C$, we denote by $C\overline{\ot} D:=(C\ot D,
\Delta_{C\ot D})$
the coalgebra in $\C$ with comultiplication
\begin{equation}\label{eq:braided-comult}
\Delta_{C\ot D}=(\id \ot \, c_{C,D} \ot \id)\circ (\Delta_C\ot \Delta_D).
\end{equation}
Let $V,W,U\in\C$ and assume $V$ is rigid; we have natural identifications
\begin{align}\label{eqn:dual-EGNO}
\begin{split}
\Hom(U\ot V,W)&\simeq \Hom(U,W\ot V^*),\\  \Hom(V^*\ot U,W)&\simeq \Hom(U,V\ot W),
\end{split}
\end{align}
see \cite[Proposition 1.10.9]{EGNO}. For instance, the top equivalence in \eqref{eqn:dual-EGNO} is:
\begin{align*}
f&\mapsto (f\ot\id_{V^*})\circ (\id_U\ot\coev_V),  &  g&\mapsto (\id_W\ot\ev_{V})\circ (g\ot \id_V)
\end{align*} 
for $f\in \Hom(U\ot V,W)$, $g\in \Hom(U,W\ot V^*)$. In particular, if $V,W\in\C$ are rigid, then for $f\in \Hom(V,W)$  
formulae \eqref{eqn:dual-EGNO} defines 
\begin{align}\label{eqn:transpose}
\,^tf\coloneqq (\ev_{W}\ot \id_{V^*})\circ (\id_{W^*}\ot f \ot \id_{V^*})\circ (\id_{W^*}\ot \coev_V) \in  \Hom(W^*,V^*)
\end{align}
These equivalences also give $(V\ot W)^*\simeq W^*\ot V^*$; with evaluation $\widetilde{\ev}_{V\ot W}$:
\begin{align}\label{eqn:dual-identification-braided}
\xymatrix{
W^*\ot V^*\ot V\ot W\ar[rrr]^{\id_{W^*}\ot \ev_V\ot \id_W\quad}&&& W^*\ot\one\ot W\simeq 
W^*\ot W \ar[r]^{\quad\qquad \qquad \ev_W}& 
\one;
}
\end{align}
compare with \eqref{eqn:dual-identification}. 

Let $(R,m,\Delta)$ be a Hopf algebra in $\C$ with $R$ rigid; then \eqref{eqn:transpose} defines morphisms 
\begin{align}\label{eqn:dual-hopf}
\,^tm\in\Hom(R^*,R^*\ot R^*) \qquad \text{and} \qquad \,^t\Delta\in\Hom(R^*\ot R^*, R^*).
\end{align} 
It follows that  $(R^*,\,^t\Delta,\,^tm)$ is again a Hopf algebra in $\C$.

\begin{definition}\label{def:braided-twists-cocycles}
Let $(R,m,\Delta)$ be a  Hopf algebra in $\C$. 

(i) A {\it braided twist} is an invertible 
morphism $\J\in\Hom(\one, R\ot R)$ that satisfies 

\begin{align}\label{eqn:twist-b}
\begin{split}
(u_R\ot \J)\star \left((\id\ot\Delta)\circ \J\right)&=(\J\ot u_R)\star\left( (\Delta\ot\id)\circ \J\right), \\ (\id\ot\eps)\circ \J&=(\eps\ot\id)\circ \J=u_R. 
\end{split}
\end{align}
When $\J\in\Hom(\one, R\ot R)$ is a braided twist, the algebra $R$ admits a  comultiplication
\begin{align}\label{eqn:br-twist-comult}
\Delta^\J:=\J\rightharpoonup\Delta \leftharpoonup \J^{-1}, 
\end{align}
and becomes a new braided Hopf algebra, denoted $R^\J$.

(ii) A {\it braided Hopf 2-cocycle} is a convolution invertible  
$\sigma\in\Hom(R\ot R,\one)$ with
\begin{align}\label{eqn:cocycle-b}
\begin{split}
(\sigma\ot\eps)\star \left(\sigma\circ(m\ot \id)\right)&=(\eps\ot\sigma)\star \left(\sigma\circ(\id\ot m)\right),\\
\sigma\circ(u_R\ot\id)&=\sigma\circ(\id\ot u_R)=\eps.
\end{split}
\end{align}
Here the first equality is in $\Hom(R^{\overline{\ot}^3},\one)$ \cf\eqref{eq:braided-comult} while in the second we identify $R\simeq R\ot \one\simeq \one\ot R$.
If $\sigma\in\Hom(R\ot R,\one)$ is a braided Hopf 2-cocycle, then the coalgebra $R$ admits an algebra structure given by
\begin{align}\label{eqn:br-cocycle-mult}
m_\sigma:=(u_R\circ\sigma)\star m\star (u_R\circ\sigma^{-1}), 
\end{align}
and becomes a new braided Hopf algebra, denoted $R_\sigma$. Indeed, \eqref{eqn:cocycle-b} assures that $m_\sigma$ is associative and it follows that $\Delta:R_\sigma\to R_\sigma\ot R_\sigma$ is still an algebra map.

As in the linear setting, we shall refer to braided Hopf 2-cocycles as braided 2-cocycles or braided  cocycles, indistinctly.
\end{definition}

Let $(R,m,\Delta)\in\C$ be a Hopf algebra. Let $\J$, resp. $\sigma$, be a braided twist, resp. Hopf 2-cocycle, on $R$. Since the multiplication, resp. the comultiplication, remains unchanged in $R^\J$, resp. $R_\sigma$, it follows that there are tensor-equivalences
\begin{align*}
&R\mod\simeq R^\J\mod, & R\comod\simeq R_\sigma\comod,
\end{align*}
compare with Theorem \ref{thm:twist-versus-tensor}. 

Assume that $R$ is rigid. If $\J\in\Hom(\one, R\ot R)$, then \eqref{eqn:transpose} defines a morphism $\,^t\J\in \Hom(R^*\ot R^*,\one)$ (identifying $\one^*\simeq \one$). Conversely,  if $\sigma\in \Hom(R\ot R,\one)$, we get  $\,^t\sigma\in\Hom(\one, R^*\ot R^*)$. The following is straightforward from the definitions.

\begin{proposition}\label{pro:dual}
Let $R\in\C$ be rigid.
\begin{enumerate} \renewcommand{\theenumi}{\roman{enumi}}\renewcommand{\labelenumi}{(\theenumi)}
\item A morphism $\sigma\in\Hom(R\ot R,\one)$ is a braided 2-cocycle on $R$ if and only if $\,^t\sigma\in\Hom(\one, R^*\ot R^*)$ is a braided twist for $R^*$.
\item  A morphism  $\J\in\Hom(\one,R\ot R)$ is a braided twist for $R$ if and only if $\,^t\J\in\Hom(R^*\ot R^*,\one)$ is a 2-cocycle on $R^*$.
\end{enumerate}
These constructions are dual to each other.
\qed
\end{proposition}

\subsection{Braided cleft objects}\label{sec:braided-cleft}
The notion of cleft object from \cf \S \ref{sec:cleft} is readily extended to the braided setting, see \cite{AF}. We fix $R\in\C$ a braided Hopf algebra and a right $R$-comodule algebra $\E$, with coaction $\rho:\E\to \E\ot R$. In turn, the coinvariants are defined as the equalizer of the diagram:
\begin{align*}
\xymatrix{\E^{\co R}\ar@{^{(}->}[rr]^{\iota} & & \E \ar@<0.4ex>@{->}[rr]^{\rho}\ar@<-0.4ex>@{->}[rr]_{\id\ot u_R}
& & \E\ot R }
\end{align*}
see \cite[Proposition 1.2]{AF};  note that equalizers exists in our category $\C$ as we assume that it is abelian. Then $\E$ is a {braided cleft extension} of $\E^{\co R}$ if there is a convolution invertible morphism of comodules $\gamma\in\Hom(R,\E)$ \cite[Definition 1.1]{AF}. 
Hence a braided cleft object is a braided cleft extension with $\E^{\co R}\simeq\one$. 
As in the standard case, there is a one-to-one correspondence between cleft objects $\E$ and braided Hopf 2-cocycles on $R$, see \cite[Theorem 1.3]{AF}. In particular, if $\gamma:\E\mapsto R$ is a section, the corresponding braided 2-cocycle on $R$ is
\begin{align}\label{eqn:sigma-gamma-braided}
\sigma=\left(m_\E\circ (\gamma\ot \gamma)\right)\star \left(\gamma^{-1}\circ m_R\right).
\end{align}

\subsection{Braided twists and braided cocycles in $\ydh$}\label{sec:braided-tc-yd}
We are mainly interested in $\C=\ydh$. We review the definitions in \S \ref{sec:braided-tc} in this context.

Let $R\in\ydh$ be a  Hopf algebra. We denote by $\cdot$, resp. $\delta$, the corresponding $H$-action, resp. $H$-coaction. Now, the product in \eqref{eq:braided-mult} reads
\begin{align}\label{eqn:braid-prod-YD}
(r\ot s)(r'\ot s')=r(s_{(-1)}\cdot r')\ot s_{(0)}s, \qquad r,r',s,s'\in R.
\end{align}
Also, recall that the vector space $R\ot H$ becomes a Hopf algebra, denoted $R\# H$, with multiplication and comultiplication given by, for $r,r'\in R, h,h'\in H$:
\begin{align}\label{eqn:bosonization}
\begin{split}
(r\# h)(r'\# h')& =r(h_{(1)}\cdot r')\ot h_{(2)}h';\\
\Delta(r\# h)&=r_{(1)}\#r_{(2)(-1)}h_{(1)}\# r_{(2)(0)}h_{(2)}.
\end{split}
\end{align}

If $R\in \ydh$, then we identify $\Hom(\k,R\ot R)\simeq R\ot R$, via $f\mapsto f(1)$.
\begin{proposition}\label{pro:dual-yd-def} Let $R\in\ydh$ be a  Hopf algebra.

(i) A braided twist for $R$ is an invertible element $\J\in R\ot R$ satisfying
\begin{align}
\label{eqn:twist-c-1} \delta(\J)=1\ot \J, \qquad h\cdot \J&=\eps_H(h)\J, \quad h\in H; 
\\
\label{eqn:twist-c-2}
\J^1\ot \mJ^1 (\mJ^2_{(-1)}\cdot \J^2_{(1)})\ot  \mJ^2_{(0)} \J^2_{(2)}&=\mJ^1(\mJ^2_{(-1)}\cdot \J^1_{(1)})\ot \mJ^2_{(0)}\J^1_{(2)}\ot \J^2; \\
\label{eqn:twist-c-3} (\eps\ot\id)(\J)&=(\id \ot\eps)(\J)=1.
\end{align}
Here $\mJ=\mJ^1\ot \mJ^2$ stands for another copy of $\J=\J^1\ot \J^2$. 

(ii) A braided 2-cocycle for $R$ is a linear map $\sigma\in \Hom(R\ot R,\k)$  which is convolution invertible and satisfies
\begin{align}\label{eqn:cocycle-c-1} 
&r_{(-1)}s_{(-1)}\sigma(r_{(0)},s_{(0)})=\sigma(r,s)u_H, \qquad \sigma(h_{(1)}\cdot r,h_{(2)}\cdot s)=\eps(h)\sigma(r,s),\\
& \label{eqn:cocycle-c-2} \sigma(x,1)=\sigma(1,x)=\eps(x) \quad \text{and }\\
\label{eqn:cocycle-c-3}
& \sigma(x_{(1)},x_{(2)(-1)}\cdot y_{(1)})\sigma(x_{(2)(0)}y_{(2)},z)
\\ &\notag \qquad \qquad \qquad \qquad \qquad \qquad  \qquad 
=\sigma(y_{(1)},y_{(2)(-1)}\cdot z_{(1)})\sigma(x,y_{(2)(0)}z_{(2)}).
\end{align}
\end{proposition}
\pf 
{\it (i)} Let $\widetilde \J\in\Hom(\k,R\ot R)$ be a braided twist for $R$, as in \eqref{eqn:twist-b} and set $\J=\widetilde \J(1)\in R\ot R$, under the identification $\Hom(\k,R\ot R)\simeq R\ot R$. Then \eqref{eqn:twist-c-1} is equivalent to the fact that $\widetilde{\J}$ is a morphism in $\ydh$. On the other hand, \eqref{eqn:twist-c-2} and \eqref{eqn:twist-c-3} are \eqref{eqn:twist-b} in this context.

{\it (ii)} Is analogous: if $\sigma\in\Hom(R\ot R,\k)$ is a braided 2-cocycle, then \eqref{eqn:cocycle-c-1} means that $\sigma$ is a morphism in $\ydh$ while \eqref{eqn:cocycle-c-2} and \eqref{eqn:cocycle-c-3} are equivalent to \eqref{eqn:cocycle-b}.
\epf
\begin{remark}
If $\J\in R\ot R$ is a braided twist, the comultiplication $\Delta^\J$ on the braided Hopf algebra $R^\J\in\ydh$ \cf \eqref{eqn:br-twist-comult} becomes 
\begin{align}\label{eqn:br-twist-comult-1}
\Delta^\J(r)=\J^1(\J^2_{(-2)}\cdot r_{(1)})(\J^2_{(-1)}r_{(2)(-1)}\cdot \J^{-1})\ot \J^2_{(0)}r_{(2)(0)}\J^{-2}.
\end{align}
Now, we can form the bosonization $R^\J\# H$ \cf \eqref{eqn:bosonization}, which is a (standard) Hopf algebra. We denote its multiplication $\Delta^\J$ by abuse of notation: by \eqref{eqn:br-twist-comult-1}:
\begin{multline}\label{eqn:DeltaJ}
\Delta^\J(r\# 1)=\J^1(\J^2_{(-2)} \cdot r_{(1)})(\J^2_{(-2)}r_{(2)(-2)} \cdot \J^{-1})
\# \J^2_{(-1)}r_{(2)(-1)}\J^{-2}_{(-1)}
\\ \ot 
\J^2_{(0)}r_{(2)(0)}\J^{-2}_{(0)}\#1, \quad r\in R.
\end{multline}
\end{remark}

We fix $R\in\ydh$ finite-dimensional. We identify $(R\ot R)^*\simeq R^*\ot R^*$ as in \eqref{eqn:dual-identification-braided}, with evaluation map $\widetilde{\ev}\coloneqq \widetilde{\ev}_{R\ot R}$ given by
\begin{align}\label{eqn:dual-ident-final}
\lg f\ot g,x\ot y \rg=\lg f,y \rg \lg g,x\rg, \quad f,g\in R^*, x,y\in R.
\end{align}
Fix a pair of (linear) dual bases $\{r_i\}_{i\in I}$ and $\{r^i\}_{i\in I}$ of $R$ and $R^*$. In particular, $\{r_i\ot r_j\}_{i,j\in I}$ and $\{r^j\ot r^i\}_{i\in I}$ is a pair of dual basis of $R\ot R$ and $R^*\ot R^*$. 
Now Proposition \ref{pro:dual} in this setting reads as follows.
\begin{proposition}\label{pro:dual-yd}
Assume $R\in\ydh$ is finite-dimensional. 
\begin{enumerate} \renewcommand{\theenumi}{\roman{enumi}}\renewcommand{\labelenumi}{(\theenumi)}
\item A morphism $\sigma\in\Hom(R\ot R,\k)$ is a braided 2-cocycle on $R$ if and only if 
\begin{align*}
\,^t\sigma(1)\in R^*\ot R^*, \qquad x\ot y\mapsto \sigma(y\ot x), \quad x,y\in R
\end{align*}
 is a braided twist for $R^*$.
\item  An element $\J\in R\ot R$ is a braided twist for $R$ if and only if 
\begin{align*}
\,^t\J\in\Hom(R^*\ot R^*,\k), \qquad f\ot g\mapsto (g\ot f)(\J), \quad f,g\in R^*   
\end{align*}
is a braided 2-cocycle for $R^*$.
\end{enumerate}
These constructions are dual to each other.\qed
\end{proposition}
\pf
This is a direct transcription of Proposition \ref{pro:dual}, by interpreting \eqref{eqn:transpose} in this context. If $\sigma\in \Hom(R\ot R,\k)$ then $\,^t\sigma\in \Hom(\k,R^*\ot R^*)$ is the composition
\begin{align}\label{eqn:dual-bases}
1&\mapsto \sum_{i,j\in I} r_i\ot r_j\ot r^j\ot r^i\mapsto \sum_{i,j\in I} \sigma(r_i\ot r_j)r^j \ot r^i,
\end{align}
so it maps $x\ot y\mapsto \sum_{i,j\in I}\sigma(r_i\ot r_j)\lg r^j,x\rg\lg r^i,y\rg=\sigma(y\ot x)$ for every $x,y\in R$.
On the one other hand, it follows from Propositions \ref{pro:dual} (i) and \ref{pro:dual-yd-def} (i) that the element $\,^t\sigma(1)\in R^*\ot R^*$ is a braided twist. Hence (i) follows. 

For (ii), we identify $\J\in R\ot R$ with $\widetilde\J\in\Hom(\k,R\ot R)$ via $\widetilde{\J}(1)=\J$. Now, \eqref{eqn:dual-identification-braided} yields a morphism $\,^t\J\coloneqq \,^t\widetilde\J\in \Hom(R^*\ot R^*,\k)$ by 
\begin{align*}
f\ot g\mapsto f\ot g \ot \J\mapsto \lg f\ot g,\J\rg=\lg g, \J^1\rg\lg f, \J^2\rg= (g\ot f)(\J), \quad f,g\in R^*.   
\end{align*}
Hence (ii) follows from Propositions \ref{pro:dual} (ii) and \ref{pro:dual-yd-def} (ii).
\epf
\subsection{The finite-dimensional setting}

Along this subsection, we fix a Hopf algebra $H$ of finite dimension. We fix, as well, a finite-dimensional Hopf algebra $R\in\ydh$ and set $A=R\# H$. Let $\J$ be a braided twist for $R\in\ydh$ and $F$ a twist for $A$.  
We show how to construct a twist $\J\# F$ for $A$ in Proposition \ref{pro:JF}. 

\medskip
Recall that we can consider $R^*$ as an object in $\ydhd$, \cf \eqref{eqn:dual}. We begin with a remark on the Hopf algebra structure of $R^*\in\ydh$ described in \eqref{eqn:dual-hopf}.

\begin{remark}\label{rem:dual-hopf}
Let $m$, resp. $\Delta$, denote the multiplication, resp. the comultiplication, of $R\in\ydh$. The bialgebra structure for $R^*\in\ydh$ in \eqref{eqn:dual-hopf} is given by, \cf \eqref{eqn:dual-ident-final}:
\begin{align}\label{eqn:RdualH}
\lg f_{(1)}\ot f_{(2)},x\ot y \rg&=f(yx), & \lg fg,x \rg&=\lg g, x_{(1)}\rg \lg f, x_{(2)}\rg,
\end{align}
for $f,g\in R^*$, $x,y\in R$. Hence the bialgebra structure for $R^*\in\ydhd$ is, see \eqref{eqn:H-to-H*}:
\begin{align}\label{eqn:RdualHdual}
\lg f_{(1)}\ot f_{(2)},x\ot y \rg&=f(xy), & \lg fg,x \rg&=\lg f, x_{(1)}\rg \lg g, x_{(2)}\rg.
\end{align}
\end{remark}

\begin{lemma}\label{lem:c} 
\begin{enumerate} \renewcommand{\theenumi}{\roman{enumi}}\renewcommand{\labelenumi}{(\theenumi)}
\item A morphism $\sigma\in\Hom(R\ot R,\k)$ is a braided 2-cocycle on $R$ if and only if 
\begin{align}\label{eqn:twist-dual}
\J(\sigma)\in R^*\ot R^*, \qquad \lg \J(\sigma),x\ot y\rg=\sigma(x,y), \quad x,y\in R
\end{align}
 is a braided twist for $R^*\in\ydhd$.
\item  An element $\J\in R\ot R$ is a braided twist for $R$ if and only if 
\begin{align}\label{eqn:cocycle-dual}
\sigma_\J\in\Hom(R^*\ot R^*,\k), \qquad f\ot g\mapsto  (f\ot g)(\J), \quad f,g\in R^*
\end{align}
is a braided 2-cocycle for $R^*\in\ydhd$.
\end{enumerate}
These constructions are dual to each other.
\end{lemma}
\pf
(i) By Proposition \ref{pro:dual-yd}, $\,^t\sigma(1)\in R^*\ot R^*$, given by $x\ot y\mapsto \sigma(y,x)$, defines a braided twist for $R^*\ot R^*$ in $\ydh$. In particular, it satisfies \eqref{eqn:twist-c-1}. By \eqref{eqn:H-to-H*}, we see that $\J(\sigma)\in R^*\ot R^*$ in \eqref{eqn:twist-dual} satisfies \eqref{eqn:twist-c-1} in $\ydhd$.
Now, on the one hand, \eqref{eqn:twist-c-3} is automatic. On the other, $\J(\sigma)$ satisfies \eqref{eqn:twist-c-2} in $\ydhd$ since $\,^t\sigma(1)$ satisfies \eqref{eqn:twist-c-2} in $\ydh$ by Remark \ref{rem:dual-hopf}.
(ii) is analogous.
\epf

\begin{corollary}\label{cor:c} 
An invertible element $\J=\J^1\ot \J^2\in R\ot R$ is a braided twist for $R\in\ydh$ if and only if $\J_{21}:=\J^2\ot \J^1$ is a braided twist for $R$ in $\ydhd$.
\end{corollary}
\pf
Combine Lemma \ref{lem:c} with Proposition \ref{pro:dual-yd}.
\epf

The next result is stated in \cite[\S 4.1]{GrM}.
\begin{lemma}\label{lem:braided-to-cocycles} 
Let $\sigma$ be a braided 2-cocycle on $R$. The assignment
\begin{align}\label{eqn:braided-to-cocycles}
x\#g\ot y\#h\longmapsto \sigma(x\ot g\cdot y)\eps(h), \quad x,y\in R, \, g,h\in H. 
\end{align}
determines an $H$-invariant Hopf 2-cocycle on $A$.
\end{lemma}
\pf 
Given the Hopf algebra structure on $A$ \cf \eqref{eqn:bosonization}, \eqref{eqn:cocycle} readily follows. Indeed, let $x,y,z\in R$, $g,h,k\in H$. Set 
${\bf x}=x\#g$, ${\bf y}=y\#h$, ${\bf z}=z\#k$. Then $\sigma({\bf x}_{(1)}, {\bf y}_{(1)})\sigma({\bf x}_{(2)}{\bf y}_{(2)}, {\bf z})$ becomes 
\begin{align*}
\sigma&(x^{(1)}\#x^{(2)}_{(-1)}g_{(1)},y^{(1)}\#y^{(2)}_{(-1)}h_{(1)})\sigma((x^{(2)}_{(0)}\#g_{(2)})(y^{(2)}_{(0)}\#h_{(2)}), {\bf z})\\
&=\sigma(x^{(1)},x^{(2)}_{(-1)}\cdot (g_{(1)}\cdot y)^{(1)})\sigma(x^{(2)}_{(0)}(g_{(1)}\cdot y)^{(2)}, g_{(2)}h\cdot z)\eps(k)\\
&\overset{\eqref{eqn:cocycle-c-3}}{=}
\sigma((g_{(1)}\cdot y)^{(1)},(g_{(1)}\cdot y)^{(2)}_{(-1)}\cdot (g_{(2)}h\cdot z)^{(1)})\sigma(x,(g_{(1)}\cdot y)^{(2)}_{(0)} (g_{(2)}h\cdot z)^{(2)})\eps(k)\\
&\overset{*}{=} \sigma(g_{(1)}\cdot y^{(1)},g_{(2)}y^{(2)}_{(-1)}h_{(1)}\cdot z^{(1)})\sigma(x,g_{(3)}\cdot y^{(2)}_{(0)} g_{(4)}h_{(2)}\cdot z^{(2)})\eps(k)\\
&= \sigma( y^{(1)},y^{(2)}_{(-1)}h_{(1)}\cdot z^{(1)})\sigma(x,(g_{(1)}\cdot y^{(2)}_{(0)})(g_{(2)}h_{(2)}\cdot z^{(2)}))\eps(k).
\end{align*}
We use the YD-condition in (*). On the other hand, $\sigma({\bf y}_{(1)}, {\bf z}_{(1)})\sigma({\bf x}, {\bf y}_{(2)}{\bf z}_{(2)})$ is
\begin{align*}
\sigma(y^{(1)}&\#y^{(2)}_{(-1)}h_{(1)},z^{(1)}\#z^{(2)}_{(-1)}k_{(1)})\sigma({\bf x}, (y^{(2)}_{(0)}\#h_{(2)})(z^{(2)}_{(0)}\#k_{(2)}))\\
&=\sigma(y^{(1)},y^{(2)}_{(-1)}h_{(1)}\cdot z^{(1)})\sigma(x,g\cdot (y^{(2)}_{(0)} h_{(2)}\cdot z^{(2)}))\eps(k).
\end{align*}
Hence \eqref{eqn:cocycle} holds using that both $\sigma$ and the multiplication in $R$ are $H$-linear.
\epf
\begin{remark}\label{rem:indep}
The proof of Lemma \ref{lem:braided-to-cocycles} does not use that $\sigma$ is $H$-colinear, \ie that $\sigma(r,s)=\sigma(r_{(0)},s_{(0)})r_{(-1)},s_{(-1)}$, \cf\eqref{eqn:cocycle-c-1}. 
\end{remark}
The next proposition shows how to construct the {\it bosonization} $\J\#F$ of a braided twist $\J$ for $R\in\ydh$ and a twist $F$ for $A$.  
\begin{proposition}\label{pro:JF}
 Let $\J$ be a braided twist for $R\in\ydh$. Then 
\begin{align}\label{eqn:def:twist-bosonization-1}
\J\#1:=\J^1\# \J^2_{(-1)}\ot \J^2_{(0)}\# 1
\end{align}
is a twist for $A=R\# H$. In particular, $R^\J\# H\simeq A^{\J\#1}$. 

More generally, if $F$ is a twist for $H$, then 
\begin{align}\label{eqn:def:twist-bosonization}
\J\#F:=FJ=F^1_{(1)}\cdot \J^1\# F^1_{(2)} \J^2_{(-1)}  \ot F^2_{(1)}\cdot \J^2_{(0)} \# F^2_{(2)}
\end{align}
is a twist for $A$, with $(R^\J\# H)^F\simeq A^{\J\# F}$.
\end{proposition}
\pf
The braided twist $\J$ determines a braided 2-cocycle for $R^*$ in $\ydhd$, via $f\ot g\mapsto \lg f\ot g, \J\rg$, $f,g\in R^*$, \cf Lemma \ref{lem:c}; so \eqref{eqn:braided-to-cocycles} defines a 2-cocycle for $A^*$:
\begin{align*}
f\#\alpha \ot g\#\beta\mapsto \lg f\ot \alpha\cdot g, \J \rg &=\lg f ,\J^1\rg \lg \alpha\cdot g, 
\J^2\rg =\lg f, \J^1\rg \lg \alpha, \J^2_{(-1)}\rg \lg g, \J^2_{(0)}\rg
\end{align*}
and thus $J=\J^1\# \J^2_{(-1)}\ot \J^2_{(0)}\# 1$ is a twist for $A$. If we write the inverse of $\J$ as $\J^{-1}\ot \J^{-2}$, then it follows that $J^{-1}=\J^{-1}\# 1\coloneqq \J^{-1}\# \J^{-2}_{(-1)}\ot \J^{-2}_{(0)}\# 1$.

To see that $R^\J\# H \simeq A^J$, we check that the (linear) identity $\id:R^\J\# H \simeq A^J$ is a coalgebra morphism, as multiplication remains unchanged on both sides, and coincides with that of $R\# A$ \cf \eqref{eqn:bosonization}. On the one hand, it is easy to see that $J\Delta(h)=\Delta(h)J$, using that $\J$ is $H$-invariant. In other words, $\Delta^J_{|H}=\Delta_{|H}$. On the other, it readily follows that $\Delta^J(r\#1)$ is given by \eqref{eqn:DeltaJ}. Indeed, 
 \begin{align*}
\Delta^J(r\#1)&=(\J\#1)\Delta(r\# 1)(\J^{-1}\#1)\\
&=(\J^1\# \J^2_{(-1)})( r_{(1)}\# r_{(2)(-1)} )(\J^{-1}\# \J^{-2}_{(-1)}) \ot \J^2_{(0)}r_{(2)(0)} \J^{-2}_{(0)}\#1 
 \end{align*}
and the claim follows using \eqref{eqn:bosonization}.

Let $F\in H\ot H$ be a twist for $H$. Now, $H\subseteq A^J$ is a Hopf subalgebra and thus $F$ becomes a twist for $A^J$ and $(A^J)^F$ is again a Hopf algebra. Observe that the comultiplication in this new Hopf algebra is determined by  $\Delta'(-)=FJ\Delta(-)(FJ)^{-1}$. We shall show that $FJ$ is a twist \ie it satisfies \eqref{eqn:twist}. In this setting this is 
\begin{align*}
(1\ot F)(1\ot J)(\id\ot\Delta)(F)(\id\ot\Delta)(J)=(F\ot 1)(J\ot 1) (\Delta\ot\id)(F)(\Delta\ot\id)(J). 
\end{align*}
Hence, it suffices to check the identities
\begin{align*}
(1\ot J)(\id\ot\Delta)(F)&=(\id\ot\Delta)(F)(1\ot J), \\
(J\ot 1) (\Delta\ot\id)(F)&=(\Delta\ot\id)(F)(J\ot 1) 
\end{align*}
which follow as $H$-invariance gives $J\Delta(h)=\Delta(h)J$, for every $h\in H$.
\epf

\section{Twists and liftings}\label{sec:twist-and-lift}
We fix throughout this section a finite-dimensional connected graded Hopf algebra $R=\oplus_{n\geq 0}R^n\in\ydh$ generated by $V=R^1$. We assume that $H$ is finite-dimensional cosemisimple and set $A=R\# H$. Hence $A$ is a quotient of $T(V)\#H$. 

Consider the canonical projection and inclusion Hopf algebra maps 
$\xymatrix{A \ar@<0.4ex>@{->}[r]^{\pi}
& H \ar@<0.4ex>@{->}[l]^{\iota}}$, \cf\eqref{eqn:bosonization}. As a result, $R=A^{\co H}=\{x\in A: (\id\ot \pi)\Delta(x)=x\ot 1\}$ and the map $\vartheta:A\to R$, $x\mapsto x_{(1)}\iota\pi\Ss(x_{(2)})$ defines a projection onto $R$. See \cite[\S 1.5]{AS} for details. We recall that the comultiplication $\Delta_R:R\to R\ot R$ can be recovered as
\begin{align}\label{eqn:comult}
\Delta_R(r)=(\vartheta\ot\id)\Delta(r\#1)=r_{(1)}\iota\pi\Ss(r_{(2)})\ot r_{(3)}.
\end{align}

Let $E\in\Cleft(A)$ be such that $E$ is also a quotient of $T(V) \# H$ 
and the section $\gamma:A\to E$ satisfies $\gamma_{|H}\in \Alg(H,E)$. This class of cleft objects naturally appears when considering lifting Hopf algebras as cocycle deformations, see \cite{AAGMV}. 

Let $\rho:E\to E\ot A$ denote the coaction.
We recall from \cite[\S 5.4]{AAGMV} that in this setting $E$ is an $H$-cleft object via $\pi$ and 
\begin{align*}
\E=E^{\co H}=\{a\in E:(\id\ot\pi)\rho(a)=a\ot 1\}\in H\mod
\end{align*} 
is an algebra such that $E=\E\#H$, identifying $H$ with its image via $T(V)\#H\twoheadrightarrow E$. This is not necessarily an algebra in $\ydh$.

The map $p:E\to E$, $a\mapsto a_{(0)}\gamma^{-1}\iota\pi(a_{(1)})$ defines a projection onto $\E$, \cite[Remark 5.6 (a)]{AAGMV}.
Moreover, as $H$ is semisimple, $\gamma(hxk)=\gamma(h)\gamma(x)\gamma(k)$ for every $h,k\in H$, $x\in A$,  and the section restricts to a morphism $\gamma_{|R}:R\to \E$ in $H\mod$ \cf \cite[Proposition 5.8 (b, c)]{AAGMV}. 

Observe that $\E\ot R$ is an algebra in $H\mod$ with multiplication
\[(a\ot x)(b\ot y)=a (x_{(-1)}\cdot b)\ot x_{(0)}y, \qquad a,b\in\E, x,y\in R.\]
Some additional properties of $\E\in H\mod$ are studied in the next lemma.

\begin{lemma}\label{lem:E}
$\E$ is an $R$-comodule algebra in $H\mod$ with $R$-coaction
\begin{align*}
\varrho(a)&=p(a_{(0)})\ot a_{(1)}=a_{(0)}\gamma^{-1}(\iota\pi(a_{(1)}))\ot a_{(2)}.
\end{align*}
The restriction $\gamma_{|R}:R\to \E$ is an $R$-colinear isomorphism.
\end{lemma}
\pf
We see that $\varrho(\E)\subseteq \E\ot R$ since $\rho(p(a))=a_{(0)}\gamma^{-1}(\iota\pi(a_{(3)}))\ot \vartheta(a_{(2)})$, $a\in E$.
Now we check that $\varrho$ is coassociative, using \eqref{eqn:comult}:
\begin{align*}
(\varrho\ot\id)\varrho(a)&=a_{(0)}\gamma^{-1}(\iota\pi(a_{(1)}))\ot a_{(2)}\Ss(\iota\pi(a_{(3)}))\ot a_{(4)}=(\id\ot\Delta_R)\varrho(a).
\end{align*}
On the other hand, we have $(\id\ot\eps)\varrho(a)=p(a)=a$, as $a\in\E$. Thus $\varrho$ defines a coaction.
The last assertion follows since $p\gamma=\gamma\vartheta$. Indeed,
\begin{align*}
\gamma\vartheta(x)=\gamma(x_{(1)}\iota\pi\Ss(x_{(2)}))=
\gamma(x_{(1)})\iota\pi\Ss(x_{(2)})=\gamma(x_{(1)})\gamma^{-1}\iota\pi(x_{(2)})
=p\gamma(x),
\end{align*}
using that $\gamma$ is $H$-invariant. Hence,
\begin{align*}
\varrho(\gamma_{|R}(x))&=p(\gamma(x\#1)_{(0)})\ot \gamma(x\#1)_{(1)}=p\gamma(x^{(1)}\#x^{(2)}_{(-1)})\ot x^{(2)}_{(0)}\#1\\
&=\gamma\vartheta(x^{(1)}\#x^{(2)}_{(-1)})\ot x^{(2)}_{(0)}\#1=\gamma\vartheta(x^{(1)})\ot x^{(2)}.
\end{align*}
Thus $\gamma_{|R}$ is $R$-colinear and the lemma follows.
\epf
If $A=R\#H$ and $E\in\Cleft(A)$ as above, then the corresponding Hopf 2-cocycle $\sigma\in Z^2(A)$, \cf\S \ref{sec:cleft}, is such that  $\sigma_{|H\ot H}=\eps\ot\eps$. We shall study this class of Hopf 2-cocycles in Proposition \ref{pro:braided-cocycle}.
%
%
%
%
Given any $\sigma\in Z^2(R\# H)$, we consider the map
\begin{align}\label{eqn:sigma_R}
 \sigma_{R}:R\ot R\to \k, \qquad x\ot y\mapsto \sigma(x\#1,y\#1), \ x,y\in R.
\end{align}
The following result, with a different formulation, is stated in \cite[\S 4.1]{GrM}; see also \cite[Lemma 4.5 \& Theorem 4.10]{ABM1}.
 We set
  \begin{multline*}
 Z^2(R)_0=\{\sigma:R\ot R\to \k \text{ convolution invertible, $H$-linear maps }\\
 \text{ that satisfy \eqref{eqn:cocycle-c-2} and \eqref{eqn:cocycle-c-3}}\}.
\end{multline*}
Hence $Z^2(R)=\{\sigma\in Z^2(R)_0:\sigma\text{ is }H\text{-colinear}\}$. The elements in $Z^2(R)_0$ have been considered in \cite[Definition 4.7]{ABM1}, with the name of {\it unital 2-cocycles}; in \loc the set $Z^2(R)_0$ is denoted $Z^2_H(R,\k)$.
\begin{proposition}\label{pro:braided-cocycle}
Let $\sigma\in Z^2(R\# H)$ be such that $\sigma_{|H\ot H}=\eps\ot\eps$. Then 
\begin{enumerate} \renewcommand{\theenumi}{\roman{enumi}}\renewcommand{\labelenumi}{(\theenumi)}
\item $\sigma\in  Z^2_H(R\# H)$.
\item $\sigma_{R}\in Z^2(R)_0$. 
\item The map $Z^2(R)_0\to Z^2_H(R\# H)$ as in \eqref{eqn:braided-to-cocycles} is an inverse to the map \eqref{eqn:sigma_R}.
\end{enumerate}
\end{proposition}
\pf
Let $E=A_{(\sigma)}\in\Cleft(A)$ be the corresponding cleft object. Then, the (linear) identity induces a section $\gamma:A\to E$ with $\gamma_{|H}\in\Alg(H,E)$; hence $\gamma$ is $H$-linear. In particular,
$\gamma(h)\gamma(x)=\gamma(hx)=\gamma(h_{(1)}x\#\Ss(h_{(2)})h_{(3)})=\gamma(h_{(1)}\cdot x\# h_{(2)})$.
Thus, using \eqref{eqn:sigma-gamma}:
\begin{align*}
\sigma(r\#t&,s\#1)=\gamma(r_{(1)}\#r_{(2)(-1)}t_{(1)})\gamma(s_{(1)}\#s_{(2)(-1)})\gamma^{-1}((r_{(2)(0)}\#t_{(2)})(s_{(2)(0)}\#1))\\
&=\gamma(r_{(1)}\#r_{(2)(-1)})t_{(1)}\gamma(s_{(1)}\#s_{(2)(-1)})\gamma^{-1}(r_{(2)(0)}t_{(2)}\cdot s_{(2)(0)}\# t_{(3)})\\
&=\gamma(r_{(1)}\#r_{(2)(-1)})\gamma(t_{(1)}\cdot s_{(1)}\# t_{(2)}s_{(2)(-1)}\Ss(t_{(4)}))\gamma^{-1}(r_{(2)(0)}t_{(3)}\cdot s_{(2)(0)}\#1)\\
&=\gamma(r_{(1)}\#r_{(2)(-1)})\gamma(t_{(1)}\cdot s_{(1)}\# (t_{(2)}\cdot s_{(2)})_{(-1)})\gamma^{-1}(r_{(2)(0)}(t_{(2)}\cdot s_{(2)})_{(0)}\#1)\\
&=\sigma(r\#1,t\cdot s\#1).
\end{align*}
Thus (i) follows.
For (ii), it is clear that $\sigma_R$ satisfies \eqref{eqn:cocycle-c-2}. Next, we observe that \eqref{eqn:cocycle} becomes \eqref{eqn:cocycle-c-3} when specialized in $x\#1$, $y\#1$, $z\#1\in A$, hence \eqref{eqn:cocycle-b} also holds. Finally, we use again that $\sigma_R$ is a convolution of $H$-linear, thus it is $H$-linear. (iii) follows from Lemma \ref{lem:braided-to-cocycles}, using Remark \ref{rem:indep}.
\epf

\begin{remark}\label{rem:liftings}
The map $\sigma_R$ in \eqref{eqn:sigma_R} is not necessarily a braided Hopf 2-cocycle in $\ydh$, that is $\sigma(r_{(-1)},s_{(-1)})r_{(0)}s_{(0)}=\sigma(r,s)$, $r,s\in R$, \cf\eqref{eqn:cocycle-c-1} does not hold in general, see Example \ref{eqn:notbraided}.

Moreover, let $V$ be a braided vector space of diagonal type with $\dim\B(V)<\infty$ and let $\Gamma$ be an abelian group with a principal YD-realization $V\in\ydga$. Every lifting of $V$ arises as cocycle deformations of $\B(V)\#H$ \cite{AGI}. However, $R=\B(V)\in\ydga$ admits no nontrivial deformations, see \cite[Theorem 6.4]{AKM}. 
\end{remark}

We now consider the element $\J(\sigma)\in R^*\ot R^*$ given by
\begin{align}\label{eqn:J-sigma}
x\ot y\longmapsto \sigma(x,y), \qquad x,y\in R.
\end{align}
Recall the structure of braided Hopf algebra of $R^*\in\ydhd$ \cf\eqref{eqn:RdualHdual}. We denote by $\delta:R^*\to H^*\ot R^*$, resp. $\delta^{(2)}:R^*\ot R^*\to H^*\ot R^*\ot R^*$, the coactions. 
\begin{corollary}\label{cor:c-1}
Let $\sigma\in Z^2(R\# H)$ be such that $\sigma_{|H\ot H}=\eps\ot\eps$. Then 
\begin{enumerate} \renewcommand{\theenumi}{\roman{enumi}}\renewcommand{\labelenumi}{(\theenumi)}
\item $\J(\sigma)\in (R^*\ot R^*)^{\co H^*}$, \ie $\delta^{(2)}(\J(\sigma))=1\ot \J(\sigma)$. 
\item $\J(\sigma)$ satisfies \eqref{eqn:twist-c-2} and \eqref{eqn:twist-c-3}.
\item $\J(\sigma)$ is a braided twist if and only if 
\begin{align}\label{eqn:trivial-action}
\alpha\cdot \J(\sigma)=\alpha(1)\J(\sigma), \qquad \alpha\in H^*.
\end{align}
\end{enumerate}
 \end{corollary}
\pf
(i) follows by \eqref{eqn:dual}, as $\sigma_R$ is $H$-linear. For (ii), we follow the lines of (the proof of) Proposition \ref{pro:dual-yd}. As in \loc we fix a pair of (linear) dual bases $\{r_i\}_{i\in I}$ and $\{r^i\}_{i\in I}$ of $R$ and $R^*$. Then, \cf\eqref{eqn:dual-bases}:
\[
\J'=\sum_{i,j\in I} \sigma(r_i\ot r_j)r^j \ot r^i\in R^*\ot R^*.
\]
satisfies \eqref{eqn:twist-c-2} and \eqref{eqn:twist-c-3} for the braided bialgebra $R^*\in\ydh$ and so $\J(\sigma)$ satisfies \eqref{eqn:twist-c-2} and \eqref{eqn:twist-c-3} for $R^*\in\ydhd$ by Remark \ref{rem:dual-hopf}. (iii) is by Definition \ref{def:braided-twists-cocycles}. 
\epf

\begin{remark}
The element $\J(\sigma)\in (R^*\ot R^*)^{\co H^*}$ does not necessarily satisfy the identity $h\cdot \J(\sigma)=\eps(h)\J(\sigma)$, $h\in H$, \cf Remark \ref{rem:liftings}.
\end{remark}

%

\subsection{The lifting context}
Next we specify our results to the setting \eqref{eq:cocycle-def} from the Introduction, to produce twists of the form $\J(\sigma)$ from Hopf 2-cocycles $\sigma$ attached to liftings of finite-dimensional Nichols algebras.

More precisely, we fix a \fd{} cosemisimple Hopf algebra $L$ and $U\in\ydL$ such that $\dim \B(U) < \infty$. We set $\H=\B(U)\# L$.
We fix $A\leq L$ a Hopf subalgebra and $W\subseteq U$ a braided subspace as in \eqref{eq:braided-subspace}, so $W\in\yda$. 

We set $R=\b(W^*)\in\ydad$, $\K=R\# A^*$. We assume \eqref{eq:cocycle-def} holds: \ie there is $\sigma\in Z^2(\K)$ such that $\L  = \K_\sigma$, so the transpose  $J(\sigma)=\,^t\sigma(1)$ of $\sigma$ is a twist for $\K^*\simeq \b(W)\# A$, hence for $\b(U)\# L$; here we identify $R^*\simeq \b(W)$. 

Assume that $\sigma_{|A^*\ot A^*}=\eps$ and consider the map $\sigma_{R}:R\ot R\to\k$ from \eqref{eqn:sigma_R}. We fix $\J(\sigma)\in R^*\ot R^*\simeq \B(W)\ot \B(W)$ as in \eqref{eqn:J-sigma}. Recall that $\J(\sigma)$ is not necessarily a braided twist, \cf Corollary \ref{cor:c-1}. We set, following \eqref{eqn:def:twist-bosonization-1}:
\[
\J(\sigma)\#1=\J(\sigma)^1\# \J(\sigma)^2_{(-1)}  \ot \J(\sigma)^2_{(0)} \# 1.
\]
\begin{proposition}\label{pro:braided}
The twist $J(\sigma)$ satisfies $J(\sigma)=\J(\sigma)\#1$.
\end{proposition}
\pf
If $x,y\in \B(W)$ and $a,b\in A$, then 
\[
\lg J(\sigma), x\#a,y\#b \rg=\sigma(x\#a,y\#b)=\sigma(x,a\cdot y)=\lg \J(\sigma),x\ot a\cdot y\rg,
\]
by the definitions of $J(\sigma)=\,^t\sigma(1)$ and $\J(\sigma)=\,^t\sigma_{R}(1)$, together with the fact that $\sigma\in Z^2_{A^*}(R\# A^*)$, by Proposition \ref{pro:braided-cocycle}. On the other hand, we have that
\[
\lg \J(\sigma),x\ot a\cdot y\rg=\lg \J(\sigma)^1, x\rg \lg \J(\sigma)^2 ,a\cdot y\rg=\lg \J(\sigma)^1, x\rg \lg \J(\sigma)^2_{(-1)},a\rg \lg \J(\sigma)^2_{(0)},y\rg
\]
as the evaluation is $A$-colinear, and thus the claim follows. 
\epf

\begin{example}\label{exa:diagonal}
Let $L$ be a Hopf algebra and let $U$ be a braided vector space with a principal realization $U\in\ydL$; let $\Gamma\leq Z(G(L))$ be as in \eqref{eqn:support}. 

We consider the following setting:
\begin{itemize}\label{page:setting}
\item[(i)] Let $\k\Gamma\subseteq A\subseteq L$ be a \fd{} cosemisimple Hopf algebra.
\item[(ii)] Let $W\subseteq U$ be a braided subspace of diagonal type such that the realization $U\in\ydL$ restricts to $W\in \yda$, with $\dim\b(W)<\infty$. 
\item[(iii)] Let $\L$ be a lifting of $W^*\in \ydad$.
    \end{itemize}
Then there is $\sigma\in Z^2(\A)$ inducing a twist $\J(\sigma)\# 1$ for $\H=\B(U)\# L$.
Indeed, $\L\simeq \A_\sigma$ is a cocycle deformation by \cite{AGI} and thus Proposition \ref{pro:braided} applies.
\begin{itemize}
\item[(iv)] Assume that \eqref{eqn:trivial-action} holds (here $H=A^*$) and let $J$ be a twist for $A$.
\end{itemize}
Then $\J(\sigma)\# J$ is a twist for $\H$, by Proposition \ref{pro:braided}. See also \S \ref{sec:diagonal-concrete}.
\end{example}


\section{Examples}\label{sec:main}

In this section we apply the ideas from \S \ref{sec:twist-and-lift} to suitable settings. 

\subsection{Twists in quantum linear spaces}\label{sec:qls}
Recall that a particularly fit setting for our approach of dualizing Hopf 2-cocycles is that of a braided vector space $V$ of diagonal type, as explained in \S \ref{sec:general}. This includes the case when $V$ is {\it quantum linear space},
 which is  the setting considered in \cite{Mo}.
We review this construction to show that it fits into the setting of Example \ref{exa:diagonal}.

Fix $\theta\in\N$, $\I=\I_\theta$. Let $(V,c)$ be a  quantum linear space with basis $(x_i)_{i\in\I}$ and braiding matrix $(q_{ij})_{i,j\in\I}$ with respect to this basis; set $N_i=\ord(q_{ii})\in\N$. Let $G$ be a finite group such that there is a principal realization $V\in\ydkg$ and fix $\Gamma\leq G$ as in \eqref{eqn:support}. 
The setting in \cite{Mo}  is:
\begin{itemize}\label{page:setting-mombelli}
\item[(i)] Let $F$ be a group such that $\Gamma\leq F\leq G$.
\item[(ii)] Let $W\subseteq V$ be an $F$-stable sub-object in $\ydkga$. 
\item[(iii)] Let $\D=(\xi_i)_{i\in\I}\cup(a_{ij})_{i\neq j\in\I}$ be an {\it $F$-invariant} family of scalars {\it compatible with $W$} \cf \cite[(3.13)--(3.15)]{Mo}.
\item[(iv)] Let $J_F=J_F^1\ot J_F^2$ be a twist for $F$.
    \end{itemize}
Out of this data, a twisting $A(V,G,W,F,J_F,\D)$ of $\b(V)\#\k G$ is built in \cite[Theorem 3.6]{Mo}. 
We obtain the following, see \S \ref{sec:mombelli-pf} for a proof.
\begin{proposition}\label{pro:J_F}
There is $\sigma\in Z^2(\b(W)\#\k\Gamma)$ such that $$A(V,G,W,F,J_F,\D)=\big(\b(V)\#\k G\big)^{\J(\sigma)\# J_F}.$$
\end{proposition}

\subsubsection{YD-pairs}
We start by constructing twists in the case $\theta=1$, which corresponds to a YD-pair $(g,\chi)$ associated to a Hopf algebra $H$. 
Let $V=\k_g^{\chi}\in\ydh$ and set $\H=\B(\k_g^{\chi})\# H$. Also, set $q=\chi(g)$, $N=\ord(q)$. 

The next proposition is well-known, see \eg \cite[Example 4.6]{Mo}, also \cite[\S 5.3.1]{GrM}.
\begin{proposition}\label{pro:J_xi}
Let $\xi\in\k$. Then  a twist for $\H$ is defined as 
\begin{align}\label{eqn:J-xi}
J_{\xi}=1\ot1+\xi\sum\limits_{k=1}^{N-1}\dfrac{1}{(N-k)_q!(k)_q!} x^kg^{N-k}\ot x^{N-k}.
\end{align}
\end{proposition}
\pf
Let $n:=\ord(g)$ and set $\Gamma=\lg g\rg\simeq C_n$, so $V$ can be realized in $\ydkga$.  
Set $\B=\B(V)\#\k\Gamma=\k[x]\#\k\Gamma/\lg x^N\rg$. Let $d\in V^*$ be such that $d(x)=1$ and $\tau\in \widehat{C_n}\simeq C_n$ be a generator. 
If $y=d\#\eps$, $h=\eps\#\tau$, then $\B^*\simeq \k[y]\#\k\widehat\Gamma/\lg y^N\rg\simeq \B$. 
Observe that $y^{r}(x^sg^t)=(r)_q!\delta_{r,s}$.

Now, any lifting of $V^*$ is of the form $\L(\xi)=\k[y]\#\k\Gamma/\lg y^{N}-\xi(1-h^N)\rg$, $\xi\in\k$. Moreover, $\L(\xi)$ is a cocycle deformation of $\B^*$ and the corresponding Hopf 2-cocycle $\sigma_\xi$ \cf \cite[\S 5.3 1.]{GrM} is given by
\begin{align}\label{eqn:sigma-xi}
\sigma_\xi(y^{i}h^k,y^{j} h^l)=\delta_{i+j,0}+\xi q^{jk}\delta_{i+j,N}
\end{align}
which coincides with the evaluation $\lg y^{i}h^k\ot y^{j} h^l,J_\xi\rg$, $J_\xi$ as in \eqref{eqn:J-xi}. Thus $J_\xi\in \B\ot \B$ is a twist for $\B$, hence for $\H$.
\epf

\begin{example}\label{eqn:notbraided}
In the notation of Proposition \ref{pro:J_xi}, let $\sigma=\sigma_\xi:\B^*\ot \B^*\to \k$ be as in \eqref{eqn:sigma-xi}, with $\xi\neq0$. Recall that $\B^*\simeq R\#\k\Gamma$, for $R=\k[y]/\lg y^N\rg$ and consider the corresponding map $\sigma_R:R\ot R\to \k$ \cf\eqref{eqn:sigma_R}. Then
\[
y_{(-1)}y^{N-1}_{(-1)}\sigma(y_{(0)}y^{N-1}_{(0)})=
g^N\sigma(y,y^{N-1})=\xi\,g^N.
\]

Hence we see that the right hand equation in \eqref{eqn:cocycle-c-1} does not hold if $N\neq n$ \ie $\sigma_R$ is not a braided cocycle in $\ydh$. An example of this situation is given if $H=\k\Gamma$, $\Gamma=C_{n}$, $n$ even. Indeed, we choose generators $g\in \Gamma$ and $\tau\in \widehat{\Gamma}\simeq \Gamma$. If $\chi=\tau^2$, then $N=\ord(q)=\frac{n}{2}$. 
We recall that $\sigma_R$ does satisfy the remaining conditions in \eqref{eqn:cocycle-c-1}--\eqref{eqn:cocycle-c-3}, 
by Proposition \ref{pro:braided-cocycle}.
\end{example}

\begin{corollary}
Assume $g^N=1$. Pick $\xi\in\k$ so that $\xi=0$ if $\chi^N\neq \eps$ and set 
\begin{align*}
\J_{\xi}=1\ot1+\xi\sum\limits_{k=1}^{N-1}\dfrac{1}{(N-k)_q!(k)_q!} x^k\ot x^{N-k}\in \B(\k_g^{\chi})\ot \B(\k_g^{\chi}).
\end{align*}
Then $\J_{\xi}$ is a braided twist for $\B(\k_g^{\chi})$ and $J_{\xi}=\J_{\xi}\#1$.
\end{corollary}
\pf
This is Proposition \ref{pro:braided}: notice that \eqref{eqn:trivial-action} holds.
\epf

\subsubsection{Quantum linear spaces}
Recall the setting (i)--(iv) from p. \pageref{page:setting-mombelli}. On the one hand, $\D$ defines a braided twist for $\b(W)$ in $\ydkga$ \cite[Theorem 3.7]{Mo}:
$$
\J_\D=\prod_{i\in\I} J_{\xi_i}\prod_{i\neq j\in\I}\exp_{q_{ij}}(a_{ij}x_i\ot x_j),
$$ 
$J_{\xi_i}$ as in \eqref{eqn:J-xi}, for $x=x_i$, and $\exp_{q_{ij}}(a_{ij}x_i\ot x_j)=\sum_{n=0}^{N_i}\frac{(a_{ij}q_{ji})^n}{(n)!_q}a_{ij}x_i^n\ot x_j^n$.
 On the other, $\k G$ is a braided Hopf algebra in $\ydkga$, with trivial 
module and comodule structure; hence $\R(V, G):=\b(V)\ot \k G$ is a braided Hopf algebra in $\ydkga$. More precisely, if $v,v'\in 
\b(V)$, $g,g'\in G$, then 
\begin{align}\label{eqn:mistake}
(v\sharp g)(v'\sharp g')=vv'\sharp  gg', \qquad \Delta(v\sharp  g)=v_{(1)}\sharp  g\ot v_{(2)}\sharp  g.
\end{align}
We denote $v\sharp g:=v\ot g\in \R(V,G)$.  Then a braided
twist for  $\R(V,G)\in\ydkga$ is
\begin{align*}
\J_{\D,F}=\J_D^1\sharp  J_F^1\ot \J_D^2\sharp  J_F^2=(\J_D^1\sharp  1\ot \J_D^2\sharp  1)(1\sharp  J_F^1\ot 1\sharp  
J_F^2),
\end{align*}
\cf \cite[Lemma 4.1]{Mo}. Set $R(V,G)=\R(V, G)\#\k\Gamma$ and, \cf \cite[(3.10)]{Mo},
\begin{align*}
\widehat{J}_{\D,F}&=\J_{\D,F}^1\# \J^2_{\D,F(-1)}\ot \J^2_{\D,F(0)}\# 1
\in  R(V,G)^{\ot2}.
\end{align*}
Now, consider the Hopf ideal 
\begin{align}\label{eqn:ideal}
 \mI=\lg 1\sharp  h\# 1- 1\sharp  1\# h|h\in\Gamma\rg\subseteq R(V,G)
\end{align}
identifying the ``two copies'' of $\Gamma$ and let $J_{\D,F}=\widehat{J}_{\D,F}+\mI$, the class of $\widehat{J}_{\D,F}$ modulo $\mI$. 
Then $J_{\D,F}$ is a twist for $\b(V)\#\k G$ and 
$$
\big(\b(V)\#\k G\big)^{J_{\D,F}}=A(V,G,W,F,J_F,\D):=\R(V, G)^{\J_{\D,F}}\#\k\Gamma/\mI.
$$ 
This is the content of  \cite[Theorem 3.6]{Mo}.
When $W=V$, $F=G$, $J_F=1_F=1\ot 1$, this twisted algebra is denoted $A(V,G,\D)$.

\begin{remark}
In \cite[\S 4, p. 430]{Mo}, \eqref{eqn:mistake} is mistyped as $(v\sharp g)(v'\sharp g')=v\, g\cdot 
v'\sharp  gg'$ but this does not define a braided Hopf algebra structure in $\b(V)\ot \k G$. This does not affect the results in \loc
\end{remark}

\smallbreak
Fix  $\sigma\in Z^2_{\Gamma}(\b(V^*)\#\k\Gamma)$. Let $J(\sigma)$, resp. $\J(\sigma)$, be the twist $J(\sigma)$ for $\b(V)\#\k\Gamma$, resp. the braided twist for $\b(V)\in\ydga$, \cf Proposition \ref{pro:braided}. Set $\J_\D:=\J_{\D,1_F}$.
\begin{lemma}\label{lem:twist-1}
There is $\sigma\in Z^2(\b(W)\#\k\Gamma)$ such that $\J_\D=\J(\sigma)$. Hence
$$A(V,G,W,F,1_F,\D)=\big(\b(V)\#\k G\big)^{J(\sigma)}.$$
In particular, the algebras $A(V,G,\D)$ can be obtained in this way.
\end{lemma}
\pf
Observe that $A(V,G,W,F,1_F,\D)\simeq A(V,G,W,\Gamma,1_F,\D)$. Indeed,
\begin{align*}
A(V,G,W,F,1_F,\D)&=\left(\B(V)\ot \k G\right)^{\J_{\D,1_F}}\#\k\Gamma/\mI\\ &\simeq 
\left(\B(V)\ot \k G\right)^{\J_{\D,1_\Gamma}}\#\k\Gamma/\mI=A(V,G,W,\Gamma,1_\Gamma,\D).
\end{align*}
Hence we may assume $F=\Gamma$. Set $A(V,G,W,\Gamma,1_F,\D)=:A$ and $B=\B(W)\#\k\Gamma$, so $B^*\simeq \b(W^*)\#\k \Gamma$. 
We have that 
${J}_{\D}:=J_{\D,1_F}=\J_D^1\# \J^2_{\D(-1)}  \ot  \J^2_{\D(0)} \# 1$ as in Proposition \ref{pro:JF} and thus the transpose of $J_\D$  determines a   Hopf 2-cocycle $\sigma\in Z^2_\Gamma(B^*)$.  Thus 
 $A=\big(\b(V)\#\k G\big)^{J(\sigma)}$ and $\J_D=\J(\sigma)$ by Proposition \ref{pro:braided}.
\epf

\subsubsection{Proof of Proposition \ref{pro:J_F}}\label{sec:mombelli-pf}
By Lemma \ref{lem:twist-1} and Proposition \ref{pro:JF} we can set:
$$A(V,G,W,F,1_F,\D)^{J_F}=\big(\b(V)\#\k G\big)^{\J(\sigma)\# J_F}.$$ The next corollary shows that this is indeed the algebra $A(V,G,W,F,J_F,\D)$. 

Set $A=A(V,G,W,F,J_F,\D)$. Recall the ideal $\mI$ in \eqref{eqn:ideal}, we identify
\begin{align}\label{eqn:ident}
\k G\simeq 1\sharp\k G\#1/\mI\subseteq A(V,G,W,F,J_F,\D), \qquad g\mapsto 1\sharp g\#1+\mI.
\end{align} 
It follows that the comultiplication in this subalgebra coincides with $\Delta^{J_F}$ and thus $(\k G)^{J_F}$ is a Hopf subalgebra. In particular $J_F^{-1}$ defines a twist for $A$. 

On the other hand, $A\simeq \big(\b(V)\#\k G\big)^{J_{\D,F}}$ and thus 
\begin{align*}
A^{J_F^{-1}}\simeq \big(\b(V)\#\k G\big)^{J_F^{-1}J_{\D,F}}\overset{(*)}\simeq \big(\b(V)\#\k G\big)^{J_{\D}},
\end{align*}
as $J_F^{-1}\widehat{J}_{\D,F}+\mI=\widehat{J}_{\D,1_F}+\mI=J_\D$,  via the identification \eqref{eqn:ident}.\qed

\medbreak

A family $\D=(\xi_i)_{i\in\I}\cup(a_{ij})_{i\neq j\in\I}$ as in p.~\pageref{page:setting} is called {\it q-symmetric} if $a_{ij}=-q_{ij}a_{ji}$, $i,j\in\I_\theta$. If $\widehat\D=\{q_{ij}a_{ji}-a_{ij}|1\leq i\neq j\leq \theta\}\cup\{\xi_i|i\in\I_\theta\}$, then $\widehat\D$ is q-symmetric.
The following corollary answers \cite[Question 4.1]{Mo}.
\begin{corollary}\label{cor:answer}
Assume $\widehat\D$ is $G$-invariant. Then $A(V,G,D)\simeq \b(V)\#\k G$.
\end{corollary}
\pf
Recall that $A=A(V,G,\D)=A(V,G,V,\Gamma,1_G,\D)$. Now, as in Lemma \ref{lem:twist-1}, the braided twist $\J_\D$ determines a Hopf 2-cocycle for $B=A^*\simeq \b(V^*)\#\k\Lambda$, $\Lambda=\widehat\Gamma$. 
Let $\{y_1,\dots,y_\theta\}$ be the dual basis of $V^*$.

There are are families of scalars $\boldsymbol\lambda=(\lambda_{ij})_{i\neq j\in \I_\theta}$ and $\boldsymbol\mu=(\mu_i)_{i\in\I_\theta}$ \cite{AS} such that $B_\sigma=\L=T(V)\#\k\Gamma/I(\boldsymbol\lambda,\boldsymbol\mu)$ where the Hopf ideal $I(\boldsymbol\lambda,\boldsymbol\mu)$ is generated by
\begin{align*}
y_iy_j-q_{ij}y_jy_i-&\lambda_{ij}(1-\chi_i\chi_j), & y_i^{N_i}-&\mu_i(1-\chi_i^{N_i}), & i\neq j\in \I_\theta.
\end{align*}
The product in $B_\sigma$ is given by
$
y_i\cdot_\sigma y_j=\sigma(y_i,y_j)(1-\chi_i\chi_j)$, $i,j\in\I_\theta$, 
while the Hopf 2-cocycle $\sigma$ is such that \cf\cite[(3.18)]{Mo}, if $1<a,b<N_i$, then for some $c_a\in\k^\times$:
\begin{align*}
\sigma(y_i,y_j)&=\lg y_i\ot y_j,\J_\D\rg=a_{ij} & 
\sigma(y_i^a,y_i^{b})&=\lg y_i^a\ot y_i^b,\J_\D\rg=
\delta_{a+b,N_i}c_a\xi_{i}.
\end{align*}
Hence it follows that $\lambda_{ij}=a_{ij}-q_{ij}a_{ji}$ and $\mu_i=\xi_ic_1$.
Thus, as $\widehat\D$ is $G$-invariant, $\lambda_{ij}\chi_i\chi_j=\lambda_{ij}$ and 
$\mu_i\chi_i^{N_i}=\mu_i$. Therefore, 
$B_\sigma\simeq B$ and thus $A\simeq \b(V)\#\k G$.
\epf

\subsection{Diagonal braidings}\label{sec:diagonal-concrete}
Let $q$ be an $N$th root of 1, some $N\geq 2$. We follow the lines of Example \ref{exa:diagonal} for a concrete braided vector space $V$ of Cartan type $A_2$ with parameter $q$. That is, $V=\big(\k\{x_1,x_2\},c^\q\big)$ with $\q= \left(\begin{smallmatrix}q& q^{-2}\\ q& \hspace*{-.4cm} q\end{smallmatrix}\right)$. Let $H$ be a cosemisimple Hopf algebra such that $V\in\ydh$ via a principal YD-realization $(\chi_i,g_i)_{i\in\I_2}$; 
$\Gamma=\lg g_1,g_2\rg\leq Z(G(H))$, set $\H=\B(V)\# H$. In particular, $V\in\ydga$ and every lifting of $V$ is a cocycle deformation $\H'_\sigma$ of $\H'=\B(V)\# \k\Gamma$ and hence it determines a twist $J(\sigma)=\J(\sigma)\# 1$ both for $\H'^*\simeq \H'$ and $\H$.

Examples of such $H$ are the group algebras $\k G_{\ell}$, for $G_\ell:=C_{\ell N}\times C_{\ell N}$, $\ell\in\N$. Also, assume $N$ satisfies $b=2N+1$, where $b:=p^n$, $p\in\N$ an odd prime, $n\in\N$. Let $\Gb=\GL_m(\F_b)\times \GL_m(\F_b)$, some $m\in\N$, so $G=G_2 \leq Z(\Gb)$ and this determines a realization $V\in\ydgb$, which restricts to $V\in\ydg$. In this setting, consider the inclusion $\Aa_m\hookrightarrow\GL_m(\F_b)$, $m\geq 4$, and set $B=\k\Aa_m$. A twist $J$ for $\Aa_m$ so that $B^J$ is not cocommutative is considered in \cite[Example 2.7]{GN}. If $\eqref{eqn:trivial-action}$ holds, then $\J(\sigma)\# J$ is a twist for $\H$.

\subsubsection{Towards an explicit formula for twists and cocycles} 

An explicit Hopf 2-cocycle is constructed in \cite[\S 5.3 1.]{GrM} for a quantum line, \cf \eqref{eqn:sigma-xi}. Another example is constructed in \cite[Proposition 5.9]{ABM1} for a quantum plane, and later generalized to quantum linear spaces in \cite[Theorem 3.9]{ABM2}. In the non-diagonal case, a Hopf 2-cocycle is defined in \cite{GM}, see \eqref{eqn:cocycle-GM} below.

Recall the section $\gamma$ involved in the presentation of the liftings \cf Lemma \ref{lem:E}. The Hopf 2-cocycles, hence the twists, can be computed explicitly if an expression of $\gamma$ is given.
We develop this idea in a particular case: when $V$ is of Cartan type $A_2$ with parameter $q=-1$ \ie $V=(\k\{ x_1,x_2\}, c^{\q})$, with $\q= \left(\begin{smallmatrix}-1& 1\\ -1&   -1\end{smallmatrix}\right)$; hence
\[
\b(V)=\k\lg x_1,x_2| x_1^2,x_2^2,x_{12}^2\rg,
\]
here $x_{12}=x_1x_2-x_2x_1$. We set $\H=\b(V)\# H$. Any lifting is a cocycle deformation of $\H$ and defines a cleft object $E=E(\bs\lambda)$, for some $\bs\lambda=(\lambda_1,\lambda_2,\lambda_{12})\in\k^3$. Set $y_{12}=y_1y_2-y_2y_1$. Then $E=\E\# H$, \cf\cite{AAG}, with
\begin{align}\label{eqn:ELambda=2}
\E=\E(\bs\lambda)=\k\lg y_1,y_2| y_1^2-\lambda_1,y_2^2-\lambda_2,y_{12}^2-\lambda_{12}\rg.
\end{align}

We write $\tx_1,\tx_2$ for another copy of the basis of $V$ and consider the canonical algebra projections $\tau:T(V)\# H\twoheadrightarrow E$, resp. $\pi:T(V)\# H\twoheadrightarrow \H$, so $x_i=\pi(\tx_i)$, $y_i=\tau(\tx_i)$, $i=1,2$. In both cases we identify $H$ with its image via the projection.

The algebras $\B(V)$ and $\E$ have PBW bases $\{x_2^ax_{12}^bx_1^c|0\leq a,b,c<2\}$ and $\{y_2^ay_{12}^by_1^c|0\leq a,b,c<2\}$, respectively. We consider the linear map $\H\to E$ defined on the corresponding bases via:
\begin{align}\label{eqn:section}
x_2^ax_{12}^bx_1^c h\longmapsto y_2^ay_{12}^by_1^ch, \qquad 0\leq a,b,c<2, \ h\in H.
\end{align} 
In particular, $\gamma_{|H}\in \Alg(H,E)$, with convolution inverse $h\mapsto \Ss(h)$.

\begin{lemma}\label{lem:section}
The $\k$-linear map $\gamma:\H\to E$ defined by \eqref{eqn:section} is a section.
\end{lemma}
\pf
It follows that $\gamma$ is a comodule morphism since the coaction $\rho:E\to E\ot \H$ is induced by the comultiplication of $T(V)\# H$ via $\rho=(\tau\ot\pi)\circ\Delta_{T(V)}$.
As $\gamma_{|H}$ is convolution invertible, then $\gamma$ is so, by \cite[Lemma 14]{T}.
\epf


Since $\sigma=\big(m\circ (\gamma\ot\gamma)\big)\star\big(\gamma^{-1}\circ m\big)$ \cf\eqref{eqn:sigma-gamma}, Lemma \ref{lem:section} allows us to compute the values $\sigma(x,y)$, $x,y\in \H$. We compute some if them in the next example.

\begin{example}\label{exa:cocycle-abelian}
It readily follows that $\gamma^{-1}(x_i)=x_ig_i^{-1}$, $i=1,2$. Also, we can see that $\gamma^{-1}(x_{12})=-(y_{12}+2y_2y_1)g_{12}^{-1}$. Hence, 	a straightforward computation yields:
\begin{align*}
	\sigma(x_i,x_j)&=\delta_{i,j}\lambda_i, \ i=1,2; & \sigma(x_{12},x_{12})&=\lambda_{12}.
\end{align*}
\end{example}


\subsection{Non-abelian examples}\label{sec:nonabelian}

We fix $n=3$ or $n=4$ and set $H=\k^{\Sym_n}$, $A=\k\Sym_n$.
Let $X=(12)^{\Sym_n}$ be  the conjugacy class of transpositions;  
$V=\k\{y_\eta:\eta\in X\}\in \ydh$ with action and coaction
\begin{align}\label{eqn:sn}
&\delta_\tau\cdot y_\eta=\delta_{\tau,\eta}y_\sigma, && y_\eta\mapsto 
\sum_{\omega\in\Sym_n}\sgn(\omega)\delta_{\omega}\ot y_{\omega^{-1}\eta\omega}, \qquad \tau\in\Sym_n, \eta\in X.
\end{align}
The corresponding Nichols algebra $\b(V)$ is denoted $\fk_n$; it is finite-dimensional. 
\begin{proposition}\label{pro:twist-nonabelian}
A twist for $\fk_n\# \k^{\Sym_n}$ is given by:
\begin{align}\label{eqn:twist-sn}
J_n=1\ot 1+\sum_{\omega\in\Sym_n}\sgn(\omega)\sum_{\eta,\tau\in X}y_\eta\delta_\omega\ot y_{\omega^{-1}\tau\omega}.
\end{align}
\end{proposition}
\pf
We have that $V^*\in\yda$ is the YD-module $\k\{x_\eta:\eta\in X_n\}$ with
\begin{align*}
&\tau\cdot x_\eta=|\tau|x_{\tau\eta\tau^{-1}}, && x_\eta\mapsto \eta\ot x_\eta, \qquad \tau\in\Sym_n, \eta\in X.
\end{align*}
Also, $\b(V^*)\simeq \fk_n$ as algebras and $\{y_\eta\}$, $\{x_\eta\}$ are dual bases of $V$, $V^*$. 
A non-trivial cocycle $\sigma\in Z^2(\fk_n\# \k{\Sym_n},\k^\times)$ is described in \cite{GM}:
\begin{align}\label{eqn:cocycle-GM}
\sigma(x_\tau\omega,x_{\tau'}\omega')&=\sgn(\omega), & \sigma(\omega,\omega')&=1, & \sigma(\cdot,\cdot)&=0\text{ elsewhere.}
\end{align}
Hence the transpose of $\sigma$ determines a twist for $\fk_n\# \k^{\Sym_n}$ given by \eqref{eqn:twist-sn}. 
\epf
Set $\J_n=1\ot 1+\sum_{\eta,\tau\in X}y_\eta \ot y_{\tau}\in\fk_n\ot \fk_n$. 
\begin{corollary}\label{cor:fkn}
Let $J_n$ be as in \eqref{eqn:twist-sn}. Then $J_n=\J_n\#1$.
\end{corollary}
\pf
This is Proposition \ref{pro:braided}.
\epf

The following example is of a different nature.

\begin{example}\label{exa:nonabelian-1}
Assume $n=4$. If $J$ is a twist for $A=\k{\Sym_4}$ (no examples if $n=3$), then 
$J=1\# J$ is also a twist for $\H=\fk_4\# \k{\Sym_4}$. 
Such a twist $J$ is considered in \cite[Remark 3.8]{GN}: it is lifted from a twist for the subgroup $C_2\times C_2\leq \Sym_4$, via the formula \eqref{eqn:twist-abelian}. Let $\A=\H^{J}$, so $A^J\subseteq \A$ is a cosemisimple Hopf subalgebra. Moreover, the authors show in \loc that $A^J$ is not cocommutative, hence neither are $\A$, $\A^*$. They also show that 
$G(A^J)\simeq \mathbb{D}_4$, so $\mathbb{D}_4\leq G(\A)$, and that it projects onto $\k\Sym_3$, hence so does $A$. In particular, $\k^{\Sym_3}\subseteq (\A^*)_0$ is a subalgebra.\qed
\end{example}

%
%

Let $\L$ be the unique (up to isomorphism) non-trivial lifting of $V\in\yda$ \cite{AG}. 

\begin{corollary}\label{cor:dual-s3}
Set $\O:=(\fk_3\# \k^{\Sym_3})^{J_3}$, then $\O\simeq \L^*$. The Hopf algebra $\O$ is not pointed, not copointed, not quasitriangular, not coquasitriangular and it has the Chevalley property.
\end{corollary}
\pf
On the one hand, $\L$ is a cocycle deformation of $\FK_3$ \cite{GMo} and the Hopf 2-cocycle is given by evaluating $J$ as in \eqref{eqn:twist-sn} \cf \cite{GM}. On the other, $\O^*$ is a pointed Hopf algebra with coradical $\Sym_3$ and $\O$ is not graded (as a Hopf algebra), hence it follows from \cite[Theorem 3.8]{AG} that $\O\simeq \L^*$. As an algebra $\O\simeq \b(V)\#\k^{\Sym_3}$, the comultiplication in $\O$ is given by $\Delta^J(-)=J\Delta(-)J^{-1}$.
The set of (non-isomorphic) simple $\O$-modules is given by 
$\{S_\eta:\eta\in\Sym_3\}$, where $S_\eta\simeq \k\{\delta_\eta\}\subseteq \O$. In particular, it follows from these formulas that 
$$S_\eta\ot S_\tau\simeq S_{\eta\tau}\not\simeq S_{\tau\eta}, \qquad \eta\neq\tau.$$ Hence $\O$ is not quasitriangular. The product of two 
simple $\O$-modules is again simple, thus $\O$ has the Chevalley property. This does not hold for $\L$-modules, as the product of 
two simple modules is not semisimple, but indecomposable. In particular, the coradical of $\O$ is not a subalgebra (plus, $\L$ has a 
simple 2-dimensional module \cite{G}). Also, $\L$ is not quasitriangular \cite{G}.
\epf

\subsubsection{Abelian extensions}

Let 
$\k^{G}\hookrightarrow H\twoheadrightarrow \mathbb \k F$ be an abelian extension. 
Assume there is $V\in\ydh$ such that the  coaction factors through $V\to \k^G\ot V\hookrightarrow H\ot V$ so $V\in 
\ydkag$. Set $\H=\b(V^*)\#\k G$ and let $\L$ be a pointed Hopf algebra over $\k G$ which is a cocycle deformation $\L=\H_\sigma$, then $\sigma$ gives rise to a twist $J(\sigma)$ for $\b(V)\# H$. See Lemma \ref{lem:matched2} below for a concrete example of this.

We recall that two subgroups $F,G$ of a group $\Sigma$ form a {\it matched pair} if the multiplication $F\times G\to 
\Sigma$ is an isomorphism. 
Every  extension $\k^{G}\hookrightarrow H\twoheadrightarrow \mathbb \k F$ arises from  a matched pair of groups $(G,F)$ and 
certain maps $\alpha:G\times F\times F\to\k$, $\beta:G\times G\times F\to\k$, so that $H\simeq \k^{G}*_{\alpha,\beta}\k F$, see \cite{M1}.

\begin{example}\label{exa:semidirect}
(a) Let $m=2,3$ and $C_m\leq\Sym_3$ be either $C_2=\lg(12)\rg$ or $C_3=\lg(123)\rg$. 
Then $(C_m,\Sym_3)$ is a matched pair of groups and $\Sigma=C_m\rtimes \Sym_3$.

\smallbreak
(b) Consider 
$C_{n+1}=\lg (12\cdots n+1)\rg\leq \Sym_{n+1}$, $\Sym_n=\{\omega|\omega(n+1)=n+1\}\leq\Sym_{n+1}$. Then $(\Sym_n,C_{n+1})$ 
is a matched pair with $C_{n+1}\Sym_n=\Sym_{n+1}$.
\end{example}

\begin{lemma}\label{lem:matched2}
Let $G=\Sym_3$ and $V\in \ydkag$ be as in \S \ref{sec:nonabelian}. 
Consider the matched pair of groups $(C_m,G)$ as in Example \ref{exa:semidirect} (a). If $H=\k^{G}*\k C_m$ is the corresponding abelian extension, then $V\in\ydh$.
\end{lemma}
\pf
Consider the 
(right) action $\dirt:\Sym_3\times C_m\to \Sym_3$, $(\eta,f)\mapsto 
f^{-1}\eta f$. Then the product in the group $\Sigma$ is given by $(f,\eta)(g,\tau)=(fg,(\eta\dirt g) \tau)$, $f,g\in C_m$, $\eta,\tau\in\k^G$.
We define an $H$-module structure on $V$ by setting $f\cdot y_\eta=y_{\eta\dirt f^{-1}}$. 
Now, we extend coaction $V\to\k^{G}\ot V$ from \eqref{eqn:sn} to $V\to H\ot V$ via the natural inclusion. Then $V\in\ydh$.
 Indeed, if $\delta$ denotes the coaction, 
\begin{multline*}
(\ad(f)\ot f\cdot \_)\circ  \delta(y_\eta)=\sum_{\omega\in\Sym_3}|\omega|f\delta_{\omega}f^{-1}\ot f\cdot y_{\omega^{-1}\eta\omega}\\
=\sum_{\omega\in\Sym_3}|\omega|\delta_{\omega\dirt f^{-1}}\ot  y_{(\omega^{-1}\eta\omega)\dirt f^{-1}}
\overset{(*)}{=}\sum_{\omega\in\Sym_3}|\omega|\delta_{\omega}\ot y_{\omega^{-1}(\eta\dirt f^{-1})\omega}=\delta(f\cdot y_\eta).
\end{multline*}
Equality (*) follows by replacing $\omega$ with $\omega\dirt f^{-1}$ as $|\omega|=|\omega\dirt f^{-1}|$.
\epf

\begin{remark}
Consider the matched pair $(\Sym_n,C_{n+1})$ from Example \ref{exa:semidirect} (b). 
There is a unique extension $H\simeq \k^{\Sym_n}* \k C_{n+1}$ \cf \cite[Theorem 8.1]{M1}. However, 
if $V$ is as in \S \ref{sec:nonabelian}, then there is no YD-structure on $V$ extending that of $\k^{\Sym_n}$. 
\end{remark}

\begin{example}\label{example3}
Let $J$ be as in \eqref{eqn:twist-sn}. By Lemma \ref{lem:matched2}, $J$ extends to a twist for $E=\fk_3\# H$ and thus a new Hopf algebra $E^J$ arises.
\end{example}

\end{document}